\documentclass{amsart}
\usepackage{hyperref}

\newtheorem{theorem}{Theorem}[section]
\newtheorem{definition}[theorem]{Definition}
\newtheorem{lemma}[theorem]{Lemma}
\newtheorem{corollary}[theorem]{Corollary}
\newtheorem{remark}[theorem]{Remark}
\newtheorem{hypothesis}[theorem]{Hypothesis {\bf H.}\hspace*{-0.6ex}}

\newcommand{\R}{{\mathbb R}}
\newcommand{\N}{{\mathbb N}}
\newcommand{\Z}{{\mathbb Z}}
\newcommand{\C}{{\mathbb C}}

\newcommand{\nn}{\nonumber}
\newcommand{\be}{\begin{equation}}
\newcommand{\ee}{\end{equation}}
\newcommand{\ul}{\underline}
\newcommand{\ol}{\overline}
\newcommand{\ti}{\tilde}

\newcommand{\spr}[2]{\langle #1 , #2 \rangle}
\newcommand{\id}{{\rm 1\hspace{-0.6ex}l}}

\newcommand{\sgn}{\mathrm{sgn}}
\newcommand{\tr}{\mathrm{tr}}
\newcommand{\im}{\mathrm{Im}}

\newcommand{\dom}{\mathfrak{D}}
\newcommand{\fdom}{\mathfrak{Q}}

\DeclareMathOperator{\Ran}{Ran}

\newcommand{\floor}[1]{\lfloor#1 \rfloor}
\newcommand{\ceil}[1]{\lceil#1 \rceil}

\newcommand{\eps}{\varepsilon}
\newcommand{\vphi}{\varphi}
\newcommand{\sig}{\sigma}
\newcommand{\lam}{\lambda}
\newcommand{\gam}{\gamma}

\numberwithin{equation}{section}

\begin{document}

\title[Relative Oscillation Theory]{Relative Oscillation Theory, Weighted Zeros of the Wronskian, and the Spectral Shift Function}

\author[H. Kr\"uger]{Helge Kr\"uger}
\address{Faculty of Mathematics\\
Nordbergstrasse 15\\ 1090 Wien\\ Austria}
\email{\href{mailto:helge.krueger@rice.edu}{helge.krueger@rice.edu}}
\urladdr{\href{http://math.rice.edu/~hk7/}{http://math.rice.edu/\~{}hk7/}}
\curraddr{Department of Mathematics, Rice University, Houston, TX 77005, USA}

\author[G. Teschl]{Gerald Teschl}
\address{Faculty of Mathematics\\
Nordbergstrasse 15\\ 1090 Wien\\ Austria\\ and International Erwin Schr\"odinger
Institute for Mathematical Physics, Boltzmanngasse 9\\ 1090 Wien\\ Austria}
\email{\href{mailto:Gerald.Teschl@univie.ac.at}{Gerald.Teschl@univie.ac.at}}
\urladdr{\href{http://www.mat.univie.ac.at/~gerald/}{http://www.mat.univie.ac.at/\~{}gerald/}}

\thanks{{\it Comm. Math. Phys. (to appear).}}
\thanks{{\it Research supported by the Austrian Science Fund (FWF) under Grant No.\ Y330}}


\begin{abstract}
We develop an analog of classical oscillation theory for Sturm--Liouville operators
which, rather than measuring the spectrum of one single operator, measures the
difference between the spectra of two different operators.

This is done by replacing zeros of solutions of one operator by weighted zeros of Wronskians of
solutions of two different operators. In particular, we show that a Sturm-type comparison theorem
still holds in this situation and demonstrate how this can be used to investigate the number
of eigenvalues in essential spectral gaps. Furthermore, the connection with Krein's
spectral shift function is established.
\end{abstract}

\maketitle

\section{Introduction}
\label{sec:int}

This paper is concerned with oscillation theory for Sturm--Liouville operators which
was originated with Sturm's celebrated memoir \cite{stu} over a hundred and seventy
years ago. Since then many extensions have been made, which are far too
numerous to be listed here. Rather we refer to the recent survey by Simon \cite{sim}
or the recent monograph by Zettl \cite{zet}.

However, there are still open problems remaining: While classical oscillation theory
works perfectly below the essential spectrum, its applicability in essential spectral
gaps is still limited. The aim of the present paper is to help filling this gap and provide a
new powerful oscillation theory which works inside essential spectral gaps.

Around 1948 Hartman and others extended classical oscillation theory to
singular Sturm--Liouville operators by showing the following
in a series of papers (\cite{har1}, \cite{har2}, \cite{har3}): For a given Sturm--Liouville
operator $H=-\frac{d^2}{dx^2}+q(x)$ on $L^2(0,\infty)$ with Dirichlet boundary
condition, let $\psi(\lam,x)$ be the solution of $H \psi = \lam \psi$ satisfying
the Dirichlet boundary condition. Then the dimension of the spectral projection $P_{(-\infty,
\lam)}(H)$ equals the number of zeros of $\psi(\lam)$ in $(0,\infty)$.
Of course a naive use of this result to count the number of eigenvalues inside
an essential spectral gap would just give $\infty-\infty$. However,
they also showed that this can be avoided by taking appropriate limits.
More precisely, the dimension of $P_{(\lam_0, \lam_1)}(H)$ can be obtained by
considering the difference of the number of zeros of the solutions $\psi(\lam_1)$ and
$\psi(\lam_0)$ inside a finite interval $(0,b)$, and performing a limit $b \to \infty$.
Further important extensions were later given by Weidmann in \cite{wd}, \cite{wd1} (see also his
lecture notes \cite{wdln}). Only recently it was shown in \cite{gst} by Gesztesy, Simon,
and one of us, that these limits can be avoided by using a renormalized version of oscillation
theory, that is, by counting zeros of Wronskians of solutions instead (see again 
\cite{sim} for a pedagogical discussion respectively \cite{kms3} for an interesting application
to minimal surfaces). That zeros of the Wronskian are related to oscillation theory is indicated
by an old paper of Leighton \cite{lei}, who noted that if two solutions have a
non-vanishing Wronskian, then their zeros must intertwine each other.

So there are several tools to count eigenvalues in essential spectral gaps and one
might ask what could be missing. To this end observe that all the previously
mentioned results have one thing in common, they all work
with only one operator $H$. One the other hand, a typical question inside an
essential spectral gap is, what kind of perturbations $H_1$ of a given operator $H_0$
produce a finite/infinite number of eigenvalues. Such kind of problems clearly
ask for an oscillation theory which compares the spectra of two given operators
$H_0$ and $H_1$. To the best of our knowledge, such an oscillation theory
was not available to date. We will demonstrate in this paper, that the right object
to look at are zeros of Wronskians of suitable solutions of $H_1$ respectively $H_0$.
In addition, since $H_1-H_0$ is not necessarily of one sign, we will weight the
zeros according to the sign of $H_1-H_0$ at a zero (see Section~\ref{sec:wz} below).
This will allow us to develop an analog of Sturm's oscillation theory comparing the
spectra of two operators including an analog of Sturm's comparison theorem, which is
applicable inside spectral gaps.

To further illustrate this, let us look at Kneser's result (\cite{kn}, see also \cite{gu} for
generalizations). Let
\begin{equation}
\tau = -\frac{d^2}{dx^2} + q(x)
\end{equation}
be a Sturm--Liouville operator on the interval $(0,\infty)$ with $q(x)\to 0$ as $x\to\infty$.
Then the bottom of the essential spectrum is an accumulation point of
eigenvalues if
\begin{equation}
\liminf_{x \to\infty} -4 x^2 q(x) >1
\end{equation}
and it is no accumulation point if
\begin{equation}
\limsup_{x \to\infty} -4 x^2 q(x) <1.
\end{equation}
The proof consists of the observation that the Euler equation
\begin{equation}
\tau_0= -\frac{d^2}{dx^2} + \frac{\mu}{x^2}
\end{equation}
is explicitly solvable plus an application of Sturm's comparison theorem between
$\tau_0$ and $\tau$.

The corresponding result within spectral gaps is Rofe-Beketov's theorem for periodic
operators (\cite{rb5}, see also the recent monograph \cite{rbk}). We state it in the version of \cite{kms},
where a detailed proof plus extensions are given. Suppose $\lim_{x \rightarrow \infty} x^2 (q(x) - q_0(x))$
exists. Then a boundary point $E_n$ of the essential spectrum is an accumulation
points of eigenvalues if 
\be
\lim_{x\to\infty} \kappa_n x^2 (q(x)-q_0(x)) >1
\ee
and no accumulation point if
\be
\lim_{x\to\infty }\kappa_n x^2 (q(x)-q_0(x)) <1,
\ee
where the critical constant $\kappa_n$ will be defined in Section~\ref{sec:app} below.

Comparing this result with Kneser's result, there is an obvious mismatch: While Kneser's
result covers all potentials which are above respectively below the
critical case $q(x)=\frac{-1}{4x^2}$, the above result only covers the
case where $q(x)$ has a precise asymptotic $q(x)=q_0(x) + \frac{c}{x^2} +o(x^{-2})$.
As a short application we will show how a comparison argument analogous to the
one for Kneser's theorem, can be used to fill this gap.

While the ingredients for such a comparison argument can essentially be found in the
results by Weidmann \cite{wd1} (see the paragraph after Theorem~\ref{thm:rosc} for a more detailed
discussion of earlier results), we still feel there is a need to advocate their use. In addition,
we hope that our novel interpretation as weighted zeros of Wronskians will lead to new relative
oscillation criteria and stimulate further research in this direction. Some possible extensions are
listed at the end of Section~\ref{sec:rot}. Furthermore, the arguments used in \cite{wd1} typically
lead to spectral estimates (compare Satz~4.1 in \cite{wd1} respectively Lemma~\ref{lem:estproj} below).
On the other hand, the results in \cite{har2} and \cite{gst} give not only estimates, but precise equalities
between spectral data and zeros of solutions respectively Wronskians of solutions. Hence, as our
main results we will establish precise equalities between spectral shifts and the number of weighted
zeros of Wronskians of certain solutions in Theorem~\ref{thm:wronskzeros} and Theorem~\ref{thmsing}.

\section{Weighted zeros of Wronskians, Pr\"ufer angles, and regular operators}
\label{sec:wz}

The key ingredient will be roughly speaking weighted zeros of Wronskians
of solutions of different Sturm--Liouville operators. However, this naive definition has a few
subtle problems which we need to clarify first. This will be done by giving an alternate
definition in terms of Pr\"ufer angles and establishing equivalence for those cases where
both definitions are well-defined.

To set the stage, we will consider Sturm--Liouville operators on $L^{2}((a,b), r\,dx)$
with $-\infty \le a<b \le \infty$ of the form
\begin{equation} \label{stli}
\tau = \frac{1}{r} \Big(- \frac{d}{dx} p \frac{d}{dx} + q \Big),
\end{equation}
where the coefficients $p,q,r$ are real-valued
satisfying
\begin{equation}
p^{-1},q,r \in L^1_{loc}(a,b), \quad p,r>0.
\end{equation}
We will use $\tau$ to describe the formal differentiation expression and
$H$ the operator given by $\tau$ with separated boundary conditions at
$a$ and/or $b$.

If $a$ (resp.\ $b$) is finite and $p^{-1},q,r$ are in addition integrable
near $a$ (resp.\ $b$), we will say $a$ (resp.\ $b$) is a \textit{regular}
endpoint.  We will say $\tau$ respectively $H$ is \textit{regular} if
both $a$ and $b$ are regular. 

For every $z\in\C\backslash\sig_{ess}(H)$ there is a unique (up to a constant) solution
$\psi_-(z,x)$ of $\tau u = z u$ which is in $L^2$ near $a$ and satisfies the
boundary condition at $a$ (if any). Similarly there is such a solution $\psi_+(z,x)$
near $b$.

One of our main objects will be the (modified) Wronskian
\be
W_x(u_0,u_1)= u_0(x)\, p(x)u_1'(x) - p(x)u_0'(x)\, u_1(x)
\ee
of two functions $u_0$, $u_1$ and its zeros. Here we think of $u_0$ and $u_1$ as two
solutions of two different Sturm--Liouville equations
\be \label{stlij}
\tau_j = \frac{1}{r} \Big(- \frac{d}{dx} p \frac{d}{dx} + q_j \Big), \qquad j=0,1.
\ee
Note that we have chosen $p_0=p_1\equiv p$ here. The case $p_0\ne p_1$ will be given
in \cite{kt2}.

Classical oscillation theory counts the zeros of solutions $u_j$, which are always
simple since $u$ and $pu'$ cannot both vanish except if $u\equiv 0$. This is
no longer true for the zeros of the Wronskian $W_x(u_0,u_1)$. In fact, it could even
happen that the Wronskian vanishes on an entire interval if $q_0=q_1$ on this interval
(cf.\ (\ref{dwr})). Such a situation will be counted as just one zero, in other words, we will
only count sign flips.

Furthermore, classical oscillation theory involves the spectrum of only one operator, but
we want to measure the difference between the spectra of two operators. That is, we
need a signed quantity. Hence we will weight the sign flips according to the sign of
$q_1-q_0$. Of course this {\em definition} is not good enough for the case
$q_1-q_0 \in L^1_{loc}$ considered here. Hence we will give a precise definition in terms
of Pr\"ufer variables next.

We begin by recalling the definition of Pr\"ufer variables $\rho_u$, $\theta_u$
of an absolutely continuous function $u$:
\be
u(x)=\rho_u(x)\sin(\theta_u(x)), \qquad
p(x) u'(x)=\rho_u(x) \cos(\theta_u(x)).
\ee
If $(u(x), p(x) u'(x))$ is never $(0,0)$ and $u, pu'$ are absolutely continuous, then $\rho_u$ is
positive and $\theta_u$ is uniquely determined once a value of
$\theta_u(x_0)$ is chosen by requiring continuity of $\theta_u$.

Notice that
\be \label{wpruefer}
W_x(u,v)= -\rho_u(x)\rho_v(x)\sin(\Delta_{v,u}(x)), \qquad
\Delta_{v,u}(x) = \theta_v(x)-\theta_u(x).
\ee
Hence the Wronskian vanishes if and only if the two Pr\"ufer angles differ by
a multiple of $\pi$. We will call the total difference
\be
\#_{(c,d)}(u_0,u_1) = \ceil{\Delta_{1,0}(d) / \pi} - \floor{\Delta_{1,0}(c) / \pi} -1
\ee
the number of weighted sign flips in $(c,d)$, where we have written $\Delta_{1,0}(x)=
\Delta_{u_1,u_0}$ for brevity.
 
Next, let us show that this agrees with our considerations above.
We take two solutions $u_j$, $j=1,2$, of $\tau_j u_j =\lam_j u_j$ and
associated Pr\"ufer variables $\rho_j$, $\theta_j$. Since we can
replace $q \to q - \lam r$ it is no restriction to assume $\lam_0=\lam_1=0$.

Under these assumptions $W_x(u_0,u_1)$ is absolutely continuous and satisfies
\be \label{dwr}
W'_x(u_0,u_1) = -(q_0(x)-q_1(x)) u_0(x) u_1(x).
\ee

\begin{lemma}\label{lem:delinc}
Abbreviate $\Delta_{1,0}(x) = \theta_1(x) - \theta_0(x)$ and suppose
$\Delta_{1,0}(x_0) \equiv 0 \mod \pi$. If $q_0(x)-q_1(x)$ is
(i) negative, (ii) zero, or (iii) positive for a.e.\ $x\in(x_0,x_0+\eps)$
respectively for a.e.\ $x\in(x_0-\eps,x_0)$ for some $\eps>0$, then the same is true for
$(\Delta_{1,0}(x) - \Delta_{1,0}(x_0))/(x - x_0)$.
\end{lemma}

\begin{proof}
By (\ref{dwr}) we have
\begin{align} \nn
W_x(u_0,u_1) &= - \rho_0(x)\rho_1(x)\sin(\Delta_{1,0}(x))\\
&=  - \int_{x_0}^x (q_0(t)-q_1(t)) u_0(t) u_1(t) dt
\end{align}
and there are two cases to distinguish: Either $u_0(x_0)$, $u_1(x_0)$ are both different
from zero or both equal to zero. If both are different from zero, $u_0(t) u_1(t)$ does not change sign
at $t=x_0$. If both are equal to zero, they must change sign at $x_0$
and hence again $u_0(t) u_1(t)$ does not change sign at $t=x_0$.  Now the claim is evident.
\end{proof}

\noindent
Hence $\#_{(c,d)}(u_0,u_1)$ counts the weighted sign flips of the Wronskian
$W_x(u_0,u_1)$, where a sign flip is counted as $+1$ if $q_0-q_1$ is positive
in a neighborhood of the sign flip, it is counted as $-1$ if $q_0-q_1$ is negative
in a neighborhood of the sign flip. If $q_0-q_1$ changes sign (i.e., it is positive
on one side and negative on the other) the Wronskian will not change its sign.
In particular, we obtain

\begin{lemma}\label{lemnbzer2}
Let $u_0$, $u_1$ solve $\tau_j u_j = 0$, $j=0,1$, where $q_0-q_1\ge 0$.
Then $\#_{(a,b)}(u_0,u_1)$ equals the number sign flips of $W(u_0,u_1)$
inside the interval $(a,b)$.
\end{lemma}

\noindent
In the case $q_0-q_1\le 0$ we get of course the corresponding negative number
except for the fact that zeros at the boundary points are counted as well since
$\floor{-x}= -\ceil{x}$. That is, if $q_0-q_1> 0$, then $\#_{(c,d)}(u_0,u_1)$
equals the number of zeros of the Wronskian in $(c,d)$ while if $q_0-q_1< 0$,
it equals minus the number of zeros in $[c,d]$. In the next theorem we will see that
this is quite natural. In addition, note that $\#(u,u) = -1$.

Finally, we establish the connection with the spectrum of regular operators.
For this let $H_0$, $H_1$ be self-adjoint extensions of $\tau_0$, $\tau_1$ with the
same separated boundary conditions.

\begin{theorem} \label{thm:reg}
Let $H_0$, $H_1$ be regular Sturm--Liouville operators associated with (\ref{stlij}) and
the same boundary conditions at $a$ and $b$. Then
\be \label{eqreg}
\dim\Ran\, P_{(-\infty, \lam_1)}(H_1) - \dim\Ran\, P_{(-\infty, \lam_0]}(H_0) =
\#_{(a,b)}(\psi_{0,\pm}(\lam_0), \psi_{1,\mp}(\lam_1)).
\ee
\end{theorem}

\noindent
The proof will be given in Section~\ref{sec:proofreg} employing interpolation between
$H_0$ and $H_1$, using $H_\eps= (1-\eps) H_0 + \eps H_1$ together with
a careful analysis of Pr\"ufer angles.

It is important to observe that in the special case $H_1=H_0$, the left-hand side
equals $\dim\Ran\, P_{(\lam_1,\lam_0)}(H_0)$ if $\lam_1>\lam_0$ and
$-\dim\Ran\, P_{[\lam_0,\lam_1]}(H_0)$ if $\lam_1<\lam_0$. This is
of course in accordance with our previous observation that
$\#(\psi_{0,\pm}(\lam_0), \psi_{1,\mp}(\lam_1))$ equals the
number of zeros in $(a,b)$ if $\lam_1>\lam_0$ while it equals minus
the numbers of zeros in $[a,b]$ if $\lam_1<\lam_0$.

\section{Relative Oscillation Theory}
\label{sec:rot}

After these preparations we are now ready to develop relative oscillation theory.
For the connections with earlier work we refer to the discussion after Theorem~\ref{thm:rosc}.

\begin{definition}\label{def:wsf}
For $\tau_0$, $\tau_1$ possibly singular Sturm--Liouville operators as in (\ref{stlij}) on $(a,b)$,
we define
\be
\underline{\#}(u_0,u_1) = \liminf_{d \uparrow b,\,c \downarrow a} \#_{(c,d)}(u_0,u_1)
\quad\mbox{and}\quad
\overline{\#}(u_0,u_1) = \limsup_{d \uparrow b,\,c \downarrow a} \#_{(c,d)}(u_0,u_1),
\ee
where $\tau_j u_j = \lam_j u_j$, $j=0,1$.

We say that $\#(u_0,u_1)$ exists, if
$\overline{\#}(u_0,u_1)=\underline{\#}(u_0,u_1)$, and write
\be
\#(u_0,u_1) = \overline{\#}(u_0,u_1)=\underline{\#}(u_0,u_1)
\ee
in this case.
\end{definition}

\noindent
By Lemma~\ref{lem:delinc} one infers that $\#(u_0,u_1)$ exists
if $q_0-\lam_0r- q_1+\lam_1r$ has the same definite sign near the endpoints $a$ and $b$.
On the other hand, note that $\#(u_0,u_1)$ might not exist
even if both $a$ and $b$ are regular, since the difference of Pr\"ufer angles might
oscillate around a multiple of $\pi$ near an endpoint. Furthermore, even if it exists,
one has $\#(u_0,u_1) = \#_{(a,b)}(u_0,u_1)$ only if there are no zeros at the
endpoints (or if  $q_0-\lam_0r-q_1+\lam_1r\ge 0$ at least near the endpoints).

\begin{remark}
Note that cases like $\underline{\#}(u_0,u_1)=-\infty$ and $\overline{\#}(u_0,u_1) = +\infty$
can occur. To construct such a situation let $\tau_j=-\frac{d^2}{dx^2} + q_j$ on $(0,\infty)$ and
$\lam_0=\lam_1=0$. Let $u_j(x)$ be the solutions satisfying a Neumann boundary condition
at $u_j'(0)=0$.

Choose $q_0(x)=0$, $q_1(x)=1$ on $(0,\pi)$ such that $\theta_1(\pi)-\theta_0(\pi)=\frac{3\pi}{2}-\frac{\pi}{2}=\pi$. Next, choose $q_0(x)=1$, $q_1(x)=0$ on $(\pi,3\pi)$ such that $\theta_1(3\pi)-\theta_0(3\pi)=\frac{3\pi}{2}-\frac{5\pi}{2}=-\pi$. Next, choose $q_0(x)=0$, $q_1(x)=1$ on $(3\pi,6\pi)$ such that $\theta_1(6\pi)-\theta_0(6\pi)=\frac{9\pi}{2}-\frac{5\pi}{2}=2\pi$. Etc.
\end{remark}

\noindent
We begin with our analog of Sturm's comparison theorem for
zeros of Wronskians. We will also establish a triangle-type inequality
which will help us to provide streamlined proofs below.
As with Sturm's comparison theorem, the proofs are elementary.

\begin{theorem}[Comparison theorem for Wronskians]\label{thm:scw}
Suppose $u_j$ satisfies $\tau_j u_j = \lam_j u_j$, $j=0,1,2$,
where $\lam_0 r -q_0 \le \lam_1 r - q_1 \le \lam_2 r - q_2$.

If $c<d$ are two zeros of $W_x(u_0,u_1)$ such that $W_x(u_0,u_1)$ does not
vanish identically, then there is at least one sign flip of $W_x(u_0,u_2)$ in $(c,d)$.
Similarly, if $c<d$ are two zeros of $W_x(u_1,u_2)$ such that $W_x(u_1,u_2)$ does not
vanish identically, then there is at least one sign flip of $W_x(u_0,u_2)$ in $(c,d)$.
\end{theorem}

\begin{proof}
Let $c$, $d$ be two consecutive zeros of $W_x(u_0,u_1)$. We first assume that $W_c(u_0,u_2)=0$
and consider $\tau_\eps = (2-\eps) \tau_1 + (\eps-1) \tau_2$, $\eps\in[1,2]$,
restricted to $(c,d)$ with boundary condition generated by the Pr\"ufer angle of $u_0$ at $c$. 
Set $u_\eps= \psi_{\eps,-}$, then we have $\Delta_{u_0,u_\eps}(c)=0$ for all $\eps$.
Moreover, we will show in (\ref{thetadot}) that $\Delta_{u_0,u_\eps}(d)$ is increasing,
implying that $W_x(u_0,u_\eps)$ has at least one sign flip in $(c,d)$ for $\eps>1$.

To finish our proof, let $\ti{u}_2$ be a second linearly independent solution. Then,
since $W(u_2,\ti{u}_2)$ is constant, we can assume $0<\Delta_{\ti{u}_2,u_2}(x)<\pi$.
This implies $\Delta_{\ti{u}_2,u_0}(c) = \Delta_{\ti{u}_2,u_2}(c)<\pi$ and
$\Delta_{\ti{u}_2,u_0}(d) = \Delta_{\ti{u}_2,u_2}(d) + \Delta_{u_2,u_0}(d)> \pi$. Consequently
$W_x(u_0,\ti{u}_2)$ also has at least one sign flip in $(c,d)$.

The second claim is proven analogously.
\end{proof}

\begin{theorem}[Triangle inequality for Wronskians]\label{thm:wtriang}
Suppose $u_j$, $j=0,1,2$ are given functions with $u_j$, $pu_j'$ absolutely continuous and
$(u_j(x), p(x) u_j'(x)) \neq (0,0)$ for all $x$. Then
\be 
\ul{\#}(u_0,u_1) + \ul{\#}(u_1,u_2) - 1 \leq \ul{\#}(u_0,u_2) \leq
\ul{\#}(u_0,u_1) + \ul{\#}(u_1,u_2) + 1
\ee
and similarly for $\ul{\#}$ replaced by $\ol{\#}$.
\end{theorem}

\begin{proof}
Take $a < c < d < b$. By definition
$$
\#_{(c,d)} (u_0, u_2) = \ceil{\Delta_{2,0} (d) / \pi} - \floor{\Delta_{2,0} (c) / \pi} - 1
$$
and using $\floor{x}+\floor{y} \le \floor{x+y} \le \floor{x}+\floor{y} +1$ respectively
$\ceil{x}+\ceil{y} -1 \le \ceil{x+y} \le \ceil{x}+\ceil{y}$ and
$\Delta_{2,0} = \Delta_{2,1} + \Delta_{1,0}$, we obtain
$$
\#_{(c,d)} (u_0, u_2) \leq \#_{(c,d)} (u_0,u_1) + \#_{(c,d)} (u_1,u_2) + 1.
$$
Thus the result follows by taking the limits $c\downarrow a$ and $d\uparrow b$.
\end{proof}

\noindent
We recall that in classical oscillation theory
$\tau$ is called oscillatory if a solution of $\tau u = 0$ has infinitely many
zeros.

\begin{definition}\label{def:relosc}
We call $\tau_1$ relatively nonoscillatory with respect to
$\tau_0$, if the quantities $\underline{\#}(u_0, u_1)$ and
$\overline{\#}(u_0, u_1)$ are finite for all solutions
$\tau_j u_j = 0$, $j = 0,1$.

We call $\tau_1$ relatively oscillatory with respect to
$\tau_0$, if one of the quantities $\underline{\#}(u_0, u_1)$ or
$\overline{\#}(u_0, u_1)$ is infinite for some solutions
$\tau_j u_j = 0$, $j = 0,1$.
\end{definition}

\noindent
Note that this definition is in fact independent of the solutions
chosen as a straightforward application of our triangle inequality
(cf.\ Theorem~\ref{thm:wtriang}) shows.

\begin{corollary}\label{cor:indisol}
Let $\tau_j u_j = \tau_j v_j = 0$, $j=0,1$. Then
\be
|\underline{\#}(u_0, u_1) - \underline{\#}(v_0, v_1)| \le 4,\quad
|\overline{\#}(u_0, u_1) - \overline{\#}(v_0, v_1)| \le 4.
\ee
\end{corollary}

\begin{proof}
By our comparison theorem we have $|\#(u_j, v_j)|\le 1$, $j =0,1$. Now use the triangle
inequality, twice.
\end{proof}

\noindent
The bounds can be improved using our comparison theorem for
Wronskians to be $\leq 2$ in the case of perturbations of definite sign.

If $\tau_0$ is nonoscillatory our definition reduces to the classical one.

\begin{lemma}
Suppose $\tau_0$ is a nonoscillatory operator, then
$\tau_1$ is relatively nonoscillatory (resp.\ oscillatory) with
respect to $\tau_0$, if and only if $\tau_1$ is
nonoscillatory (resp.\ oscillatory).
\end{lemma}

\begin{proof}
This follows by taking limits in
$$
|\#_{(c,d)}(u,v) - \#_{(c,d)}(u) + \#_{(c,d)}(v)|\leq 2,
$$
where $\#_{(c,d)}(u)= \ceil{\theta_u(d) / \pi} - \floor{\theta_u(c) / \pi} -1$ is the
number of zeros of $u$ inside $(c,d)$.
\end{proof}

\noindent
To demonstrate the usefulness of Definition~\ref{def:relosc},
we now establish its connection with the spectra
of self-adjoint operators associated with $\tau_j$, $j=0,1$.

\begin{theorem} \label{thm:rosc}
Let $H_j$ be self-adjoint operators associated with $\tau_j$, $j=0,1$. Then
\begin{enumerate}
\item
$\tau_0-\lam_0$ is relatively nonoscillatory with respect to $\tau_0-\lam_1$
if and only if $\dim\Ran P_{(\lam_0,\lam_1)}(H_0)<\infty$.
\item
Suppose $\dim\Ran P_{(\lam_0,\lam_1)}(H_0)<\infty$ and
$\tau_1-\lam$ is relatively nonoscillatory with respect to $\tau_0-\lam$ for one
$\lam \in [\lam_0,\lam_1]$.  Then it is relatively nonoscillatory for
all $\lam\in[\lam_0,\lam_1]$ if and only if $\dim\Ran P_{(\lam_0,\lam_1)}(H_1)<\infty$.
\end{enumerate}
\end{theorem}

\begin{proof}
(i) This is item (i) of \cite[Thm.7.5]{gst}.
(ii) Let $\lam, \ti{\lam} \in [\lam_0,\lam_1]$,  $\tau_j u_j(\lam) = \lam u_j(\lam)$, $j=0,1$
and suppose we are relatively nonoscillatory at $\lam$.
Then applying our triangle inequality twice implies
\begin{align*}
\overline{\#}(u_0(\ti{\lam}),u_1(\ti{\lam})) &\leq
\overline{\#}(u_0(\ti{\lam}), u_1(\lam)) + \#(u_1(\lam),u_1(\ti{\lam})) +1\\
&\leq  \#(u_0(\ti{\lam}), u_0(\lam)) + \overline{\#}(u_0(\lam),u_1(\lam)) +
\#(u_1(\lam),u_1(\ti{\lam})) + 2,
\end{align*}
and similar estimates with the roles of $\lam$ and $\ti{\lam}$
interchanged and $\overline{\#}$ replaced by $\underline{\#}$.
Hence  if $\dim\Ran P_{(\lam_0,\lam_1)}(H_1)<\infty$ we are relatively
nonoscillatory by (i). The converse direction is proven analogously.
\end{proof}

\noindent
We remark that item (i), which corresponds to the case of equal operators $H_0=H_1$
but different spectral parameters $\lam_0\ne \lam_1$, is what was called
renormalized oscillation theory by Gesztesy, Simon, and Teschl in \cite{gst}. It
also follows from earlier results by Hartman \cite{har2}, Weidmann \cite[Satz~4.1]{wd1}
(see also \cite{wdln}), which was pointed out in the appendix of \cite{kms} and was called
relative oscillation theorem there. An argument in the spirit of item (ii) was also one of the
key ingredients in Rofe-Beketov's original work (cf.\ \cite{rb5}).

For a practical application of this theorem one needs of course criteria
when  $\tau_1-\lam$ is relatively nonoscillatory with respect to $\tau_0-\lam$
for $\lam$ inside an essential spectral gap.

\begin{lemma}\label{lem:nonoscingap}
Let $\lim_{x\to a} r(x)^{-1} (q_0(x) - q_1(x)) = 0$ if $a$ is singular, and
similarly, $\lim_{x\to b} r(x)^{-1} (q_0(x) - q_1(x)) = 0$ if $b$ is singular.
Then $\sig_{ess}(H_0)=\sig_{ess}(H_1)$ and
$\tau_1 - \lam$ is relatively nonoscillatory with respect
to $\tau_0 - \lam$ for $\lam \in \R \backslash\sigma_{ess}(H_0)$.
\end{lemma}

\begin{proof}
Since $\tau_1$ can be written as $\tau_1 = \tau_0 + \ti{q}_0 + \ti{q}_1$, where
$\ti{q}_0$ has compact support near singular endpoints and $|\ti{q}_1|<\eps$,
for arbitrarily small $\eps>0$, we infer that $R_{H_1}(z) - R_{H_0}(z)$ is the norm
limit of compact operators. Thus $R_{H_1}(z) - R_{H_0}(z)$ is compact and hence
$\sig_{ess}(H_0)=\sig_{ess}(H_1)$.

Let $\delta>0$ be the distance of $\lam$ to the essential spectrum and choose
$a < c < d < b$, such that
$$
|r^{-1}(q_1(x) - q_0(x))| \le \delta/2,\qquad x\not\in(c,d).
$$
Clearly $\#_{(c,d)} (u_0,u_1) < \infty$, since both operators are regular
on $(c,d)$. Moreover, observe that
$$
q_0 - \lam_+ r \leq q_1- \lam r \leq q_0- \lam_- r, \qquad
\lam_\pm=\lam\pm\delta/2,
$$
on $I = (a,c)$ or $I =(d,b)$. Then Theorem~\ref{thm:rosc}~(i) implies
$\#_I(u_0(\lam_-), u_0(\lam_+)) < \infty$ and invoking Theorem~\ref{thm:scw}
shows $\#_I(u_0(\lam_\pm), u_1(\lam)) < \infty$.
From Theorem~\ref{thm:wtriang} and \ref{thm:rosc}~(i) we infer
$$
\ol{\#}_I(u_0(\lam),u_1(\lam)) < \#(u_0(\lam),u_0(\lam_+)) +
\#_I(u_0(\lam_+), u_1(\lam)) + 1 < \infty
$$
and similarly for $\ul{\#}_I(u_0(\lam),u_1(\lam))$.
This shows that $\tau_1-\lam$ is relatively nonoscillatory with respect to $\tau_0$.
\end{proof}

\noindent
Our next task is to reveal the precise relation between the number of
weighted sign flips and the spectra of $H_1$ and $H_0$. The special
case $H_0=H_1$ is covered by \cite{gst}:

\begin{theorem}[\cite{gst}] \label{thm:gst}
Let $H_0$ be a self-adjoint operator associated with $\tau_0$ and suppose
$[\lam_0,\lam_1]\cap\sig_{ess}(H_0)=\emptyset$.
Then
\be
 \dim\Ran P_{(\lam_0, \lam_1)} (H_0) = \#(\psi_{0,\mp}(\lam_0),\psi_{0,\pm}(\lam_1)).
\ee
\end{theorem}

\noindent
We will provide an alternate proof in Section~\ref{sec:appreg}.

Combining this result with our triangle inequality already gives some rough
estimates in the spirit of Weidmann \cite{wd1} who treats the case $H_0=H_1$.

\begin{lemma}\label{lem:estproj}
For $j=0,1$ let $H_j$ be a self-adjoint operator associated with $\tau_j$ and separated
boundary conditions. Suppose 
that $(\lam_0, \lam_1)\subseteq\R\backslash(\sig_{ess}(H_0)\cup\sig_{ess}(H_1))$, then
\begin{align}
\nonumber
\dim\Ran P_{(\lam_0, \lam_1)} (H_1) &-
\dim\Ran P_{(\lam_0, \lam_1)} (H_0)\\
& \leq \underline{\#}(\psi_{1,\mp} (\lam_1), \psi_{0,\pm} (\lam_1)) - 
\overline{\#}(\psi_{1,\mp} (\lam_0), \psi_{0,\pm} (\lam_0)) + 2,
\end{align}
respectively,
\begin{align}
\nonumber
\dim\Ran P_{(\lam_0, \lam_1)} (H_1) &-
\dim\Ran P_{(\lam_0, \lam_1)} (H_0)\\
& \geq \overline{\#}(\psi_{1,\mp} (\lam_1), \psi_{0,\pm} (\lam_1)) - 
\underline{\#}(\psi_{1,\mp} (\lam_0), \psi_{0,\pm} (\lam_0)) - 2.
\end{align}
\end{lemma}

\begin{proof}
By the triangle inequality (cf.\ Theorem~\ref{thm:wtriang}) we have
\begin{align*}
& \#_{(c,d)}(\psi_{1,-}(\lam_1),\psi_{1,+}(\lam_0)) - \#_{(c,d)}(\psi_{0,-}(\lam_1),\psi_{0,+}(\lam_0))\\
&\qquad \le \#_{(c,d)}(\psi_{1,-}(\lam_1),\psi_{0,+}(\lam_1)) + \#_{(c,d)}(\psi_{1,-}(\lam_1),\psi_{0,+}(\lam_1)) +2.
\end{align*}
The result now follows by taking limits using that
$$
\lim_{c\downarrow a, d \uparrow b} \#_{(c,d)}(\psi_{1,-}(\lam_1),\psi_{1,+}(\lam_0))
= \dim\Ran P_{(\lam_0, \lam_1)} (H_1)
$$
and
$$
\lim_{c\downarrow a, d \uparrow b} \#_{(c,d)}(\psi_{0,-}(\lam_0),\psi_{0,+}(\lam_1))
= -\dim\Ran P_{(\lam_0, \lam_1)} (H_0)
$$
by the previous theorem. The second claim follows similarly.
\end{proof}

\noindent
To turn the inequalities into equalities in Lemma~\ref{lem:estproj}
will be one of our remaining goals.

Observe that for semibounded operators can choose $\lam_0$ below the spectra
of $H_0$ and $H_1$, causing $\#(\psi_{1,\mp}(\lam_0),\psi_{0,\pm}(\lam_0))$
to vanish:

\begin{lemma}
Let $r^{-1}(q_0- q_1)\leq \delta$ for some $\delta\in\R$. Furthermore, suppose
that the operator $H_0$ associated with $\tau_0$ is bounded from below, $H_0 \ge E_0$,
and the form domains of $H_0$ and $H_1$ are equal. Then 
\be
\#(\psi_{1,\mp}(\lam),\psi_{0,\pm}(\lam)) = 0, \quad\lam<E_0-\delta.
\ee
\end{lemma}

\begin{proof}
Suppose there is $c\in (a,b)$ such that $W_c(\psi_{0,+}(\lam),\psi_{1,-}(\lam)) = 0$.
Then, there is a $\gam$ such that
$$
\vphi(x) = \begin{cases} \psi_{1,-}(\lam,x), & x\leq c,\\
\gam\, \psi_{0,+}(\lam,x), & x \geq c, \end{cases}
$$
is continuous and hence in the form domain of $H_0$ (see the remark after Theorem~A.3 in
\cite{gst}). But then
$$
\spr{\vphi}{H_0\vphi} \leq (\lam + \delta) \|\vphi\|^2 < E_0 \|\vphi\|^2,
$$
contradicting $H_0\ge E_0$.
\end{proof}

\noindent
Our first approach will use approximation by regular problems. However, the
standard approximation technique only implies strong convergence, which
is not sufficient for our purpose. Hence our argument is based on a refinement of a
method by Stolz and Weidmann \cite{sw} which will provide convergence
of spectral projections in the trace norm for suitably chosen regular operators 
(see \cite{wd2} for a nice overview).

\begin{theorem}\label{thm:wronskzeros}
Let $H_0$, $H_1$ be self-adjoint operators associated with $\tau_0$, $\tau_1$,
respectively, and separated boundary conditions. Suppose 
\begin{enumerate}
\item
$q_1 \le q_0$, near singular endpoints,
\item
$\lim_{x\to a} r(x)^{-1}(q_0(x)-q_1(x)) = 0$ if $a$ is singular and
$\lim_{x\to b} r(x)^{-1}(q_0(x)-q_1(x)) = 0$ if $b$ is singular,
\item
$H_0$ and $H_1$ are associated with the same boundary conditions near $a$ and $b$,
that is, $\psi_{0,-}(\lam)$ satisfies the boundary condition of $H_1$ at $a$ (if any) and
$\psi_{1,+}(\lam)$ satisfies the boundary condition of $H_0$ at $b$ (if any).
\end{enumerate}

Suppose $\lam_0 < \inf\sig_{ess}(H_0)$. Then
\be\label{eq:wronskzeros0}
\dim\Ran P_{(-\infty,\lam_0)}(H_1) - \dim\Ran P_{(-\infty,\lam_0]}(H_0)
= \#(\psi_{1,\mp} (\lam_0), \psi_{0,\pm} (\lam_0)).
\ee

Suppose $\sig_{ess}(H_0) \cap [\lam_0,\lam_1] = \emptyset$. Then $\tau_1-\lam_0$ is
nonoscillatory with respect to $\tau_0-\lam_0$ and
\begin{align}
\nonumber
&\dim\Ran P_{[\lam_0, \lam_1)} (H_1) - \dim\Ran P_{(\lam_0, \lam_1]} (H_0)\\ \label{eq:wronskzeros1}
& \qquad = \#(\psi_{1,\mp} (\lam_1), \psi_{0,\pm} (\lam_1)) - 
\#(\psi_{1,\mp} (\lam_0), \psi_{0,\pm} (\lam_0)).
\end{align}
\end{theorem}

\noindent
The proof will be given in Section~\ref{sec:appreg}.

\begin{remark}
Note that condition (ii) implies $\sig_{ess}(H_0)=\sig_{ess}(H_1)$ as pointed out in
Lemma~\ref{lem:nonoscingap}.
In addition, (ii) implies that any function which is in $\dom(\tau_0)$ near $a$ (or $b$)
is also in $\dom(\tau_1)$ near $a$ (or $b$), and vice versa. Hence condition (iii) is
well-posed.
\end{remark}

\noindent
Our second approach will connect our theory with Krein's spectral
shift function. Given the regular case in Theorem~\ref{thm:reg}, we can extend this result
to operators whose resolvent difference is trace class by replacing the left-hand side in
(\ref{eqreg}) by the spectral shift function. In order to fix the unknown constant in the spectral
shift function, we will require that $H_0$ and $H_1$ are connected via a path within the set of
operators  whose resolvent difference with $H_0$ are trace class. Hence we will require

\begin{hypothesis} \label{hyp:h0h1}
Suppose $H_0$  and $H_1$ are self-adjoint operators associated with $\tau_0$ and $\tau_1$
and separated boundary conditions. Abbreviate $\ti{q} = r^{-1}|q_0 - q_1|$, and assume that:
\begin{itemize}
\item[(i)] $\ti{q}$ is relatively bounded with respect to $H_0$ with $H_0$-bound less than one or
\item[(i')] $H_0$ is bounded from below and $\ti{q}$ is relatively form bounded with respect
to $H_0$ with relative form bound less than one and
\item[(ii)] $\sqrt{\ti{q}} R_{H_0} (z)$ is Hilbert--Schmidt for one (and hence for all)
$z\in\rho(H_0)$.
\end{itemize}
\end{hypothesis}

\noindent
It will be shown in Section~\ref{sec:xi} that these conditions ensure that we can interpolate
between $H_0$ and $H_1$ using operators $H_\eps$, $\eps\in[0,1]$, such that the resolvent difference
of $H_0$ and $H_\eps$ is continuous in $\eps$ with respect to the trace norm. Hence we can fix
$\xi(\lam,H_1,H_0)$ by requiring $\eps \mapsto \xi(\lam,H_\eps,H_0)$ to be continuous in
$L^1(\R,(\lam^2+1)^{-1}d\lam)$, where we of course set $\xi(\lam,H_0,H_0)=0$
(see Lemma~\ref{lem:resolvconv}). While $\xi$ is only defined a.e., it is constant on the intersection
of the resolvent sets $\R\cap\rho(H_0)\cap\rho(H_1)$, and we will require it to be continuous there.
In particular, note that by Weyl's theorem the essential spectra of $H_0$ and
$H_1$ are equal, $\sig_{ess}(H_0)=\sig_{ess}(H_1)$. Then we have the following result:

\begin{theorem} \label{thmsing}
Let $H_0$, $H_1$ satisfy Hypothesis~\ref{hyp:h0h1}. Then
for every $\lam\in\R\cap\rho(H_0)\cap\rho(H_1)$ we have
\be
\xi(\lam,H_1,H_0) =
\#(\psi_{0,\pm}(\lam), \psi_{1,\mp}(\lam)).
\ee
\end{theorem}

\noindent
The proof will be given in Section~\ref{sec:apptr}.

In particular, this result implies that under these assumptions $\tau_1-\lam$ is relatively nonoscillatory
with respect to $\tau_0-\lam$ for every $\lam$ in an essential spectral gap.

The main idea is to interpolate between $H_0$ and $H_1$.
Under proper assumptions it is then possible to control $\xi(\lam,H_\eps,H_0)$.
However, it seems extremely hard to control the zeros of the Wronskian. To do this we will have to
assume that $q_1-q_0$ has compact support. We will then remove this restriction by extending the
support first to one and then to the other side. The details will be given in Section~\ref{sec:apptr}.

Finally, we remark that since the results from \cite{gst} extend to one-dimensional
Dirac operators \cite{toscd} (see also \cite{kms2}) and Jacobi operators
\cite{tosc} (compare also \cite[Chap.~4]{tjac}), similar results are expected to
hold for these operators and will be given in \cite{troscd} respectively \cite{troscj}.
Furthermore, it will of course be interesting to develop relative oscillation criteria
for Sturm--Liouville operators! What are the analogs of some classical oscillation criteria?
Is there an analog of \cite{gu}? Some results concerning these questions will be given in
\cite{kt3}.

\section{Applications}
\label{sec:app}

In this section we want to look at the classical problem of the number of eigenvalues
in essential spectral gaps of perturbed periodic operators \cite{bi}, \cite{rb1} (see also \cite{gs}),
\cite{zhe}. The precise critical case was first determined by Rofe-Beketov in a series of papers
\cite{rb2}--\cite{rb5} with later additions by Khryashchev \cite{khr} and Schmidt \cite{kms}.
The purpose of this section is to show how the the results in \cite{kms} can be extended using
our methods.

For convenience of the reader, we recall some basic facts of the theory
of periodic differential operators first (see for example \cite{ea} or \cite{wdln}).

Let $p$, $q_0$ be $\alpha$-periodic, that is, $p(x+\alpha)=p(x)$, $q_0(x+\alpha)=q_0(x)$,
and consider the corresponding Sturm--Liouville
expressions
\be
\tau_0 = - \frac d{dx} p \frac d{dx} + q_0.
\ee
Denote by $c(\lam,x)$, $s(\lam,x)$ a fundamental system of solutions
corresponding to the initial conditions $c(\lam,0)=p(0) s'(\lam,0)=1$, $s(\lam,0)= p(0) c'(\lam,0)=0$.
In particular, their Wronskian reads $W(c,s) = 1$. We then call
\be
M(\lam) =
\begin{pmatrix} c(\lam,\alpha) & s(\lam,\alpha) \\
p(0) c'(\lam,\alpha) & p(0) s'(\lam,\alpha) \end{pmatrix}
\ee
the monodromy matrix.
The discriminant $D(\lam)$ is given by $D(\lam) = \tr(M(\lambda))$.

Since we are only interested in the question wether the number of eigenvalues
are finite or not, it suffices to look at the half-line case $(1,\infty)$ with
a, for example, Dirichlet boundary condition at $1$. Denote the corresponding self-adjoint
operator by $H_0$. The essential spectrum of $H_0$ is given by
\be
\sigma_{ess}(H_0) = \sig_{ac}(H_0) = \{ \lam\,|\, |D(\lam)| \leq 2 \}
= \bigcup_{n=0}^\infty [E_{2n},E_{2n+1}].
\ee
The critical coupling constant at the endpoint $E_n$ of an essential spectral gap
introduced in Section~\ref{sec:int} is given by (\cite{kms})
\be \label{defkap}
\kappa_n = \frac{\alpha^2}{4 |D|'(E_n)}.
\ee
It is related to the effective mass $m(E_n)$ used in solid state physics via
$\kappa_n = (8 m(E_n))^{-1}$ (see \cite{rb5}).

We note that $\kappa_{2n+1}> 0$ for a lower endpoint of a spectral gap, and $\kappa_{2n} < 0$
for an upper endpoint.

Before coming to our applications, we ensure that our hypothesis from the previous section are
satisfied. We begin by computing the form domain of $H_0$.

\begin{lemma}
Abbreviate
\be
\fdom= \{ f \in L^2(1,\infty) \,|\, f \in AC[1,\infty), \: \sqrt{p} f' \in L^2(1,\infty)\}.
\ee
The form domain of $H_0$ is given by
\be
\fdom(H_0)= \{ f \in \fdom \,|\, f(1)=0\}.
\ee
Moreover, for every $\eps>0$ there is some $C>0$ such that
\be \label{qdperbd}
\sup_{x_0 \le x \le x_0+\alpha} |f(x)|^2 \le \eps \int_{x_0}^{x_0+\alpha} p(x) |f'(x)|^2 dx
+ C  \int_{x_0}^{x_0+\alpha} |f(x)|^2 dx, \quad f\in \fdom.
\ee
In particular, $q_1-q_0$ is infinitesimally form bounded with respect to $H_0$ if
\be
\int_{x_0}^{x_0+\alpha} |q_1(x) - q_0(x)| dx < C_0, \qquad x_0 \in (1,\infty),
\ee
where $C_0$ is independent of $x_0$.
\end{lemma}

\begin{proof}
Equation (\ref{qdperbd}) is a standard Sobolev estimate. For the case
$p(x)\ne 1$ required here compare for example \cite[Lem~A.2]{gst}.

Next, set
\be
A = \sqrt{p} \frac{d}{dx}, \qquad \dom(A) = \{ f \in \fdom \,|\, f(1)=0\}
\ee
and note that $A$ is then a closed operator with adjoint given by
\begin{align}\nn
A^* &= -\frac{d}{dx} \sqrt{p} ,\\
\dom(A^*)  &= \{ f \in L^2(1,\infty) \,|\, \sqrt{p} f \in AC[1,\infty), \: (\sqrt{p} f)' \in L^2(1,\infty)\}.
\end{align}
Hence,
\begin{align}\nn
A^* A &= -\frac{d}{dx} p \frac{d}{dx},\\
\dom(A^* A) &= \{ f \in L^2(1,\infty) \,|\, f,pf' \in AC[1,\infty), \: (p f')' \in L^2(1,\infty), \: f(1)=0\}
\end{align}
is self-adjoint with $\fdom(A^* A)= \dom(A)$ and by (\ref{qdperbd})
$q_0$ is infinitesimally form bounded with respect to $A^* A$.
Since the same is true for $q_1$ by assumption, the lemma is proven.
\end{proof}

\begin{lemma}\label{lem:HS}
Let $H_0$ be an arbitrary Sturm--Liouville operator on $(a,b)$.
Then $\sqrt{q} R_{H_0}(z)$ is Hilbert--Schmidt if and only if
\be
\| \sqrt{q} R_{H_0}(z) \|_{\mathcal{J}_2}^2 = \frac{1}{\im(z)}
\int_a^b |q(x)| \im(G_0(z,x,x)) r(x) dx
\ee
is finite. Here $G_0(z,x,y) =(H_0-z)^{-1}(x,y)$ denotes the Green's function of $H_0$.
\end{lemma}

\begin{proof}
From the first resolvent identity we have
$$
G_0(z,x,y) - G_0(z',x,y) = (z-z') \int_a^b G_0(z,x,t) G_0(z',t,y) r(t) dt.
$$
Setting $x=y$ and $z'=z^*$ we obtain
$$
\im(G_0(z,x,x)) = \im(z) \int_a^b |G_0(z,x,t)|^2 r(t) dt.
$$
Using this last formula to compute the Hilbert--Schmidt norm proves the lemma.
\end{proof}

We recall that in the case of periodic operators $G(z,x,x)$ is a
bounded function of $x$. In fact, we have
\be\label{estgfper}
|G(z,x,y)| \leq const(z) \exp(-\gam(z) |x-y|),
\ee
where $\gam(z)>0$ denotes the Floquet exponent.

Now we are ready to apply our theory:

\begin{theorem}\label{thm:app}
Let $p, q_0$ $\alpha$-periodic and $q_1$ a perturbed
potential which is regular at $1$, such that either $q_0 - q_1 \in L^1(1,\infty)$ or
$\lim_{x\to\infty} (q_0(x) - q_1(x))=0$.
Define the differential expressions on $(1, \infty)$ by
\be
\tau_0 = - \frac d{dx} p(x) \frac d{dx} + q_0(x), \quad
\tau_1 =  - \frac d{dx} p(x) \frac d{dx} + q_1(x)
\ee
and let $H_0$, $H_1$ be the corresponding self-adjoint operators $L^{2}(1,\infty)$.

Let $E_n$ be an endpoint of a gap in the essential spectrum of $H_0$ with corresponding
$\kappa_n$ given by (\ref{defkap}).
Then $E_n$ is an accumulation point of eigenvalues of $H_1$ if
\be
\liminf_{x \to\infty} \kappa_n x^2 (q_1 (x) - q_0(x)) > 1.
\ee
and $E_n$ is no accumulation point of eigenvalues of $H_1$ if
\be
\limsup_{x \to\infty} \kappa_n x^2 (q_1 (x) - q_0(x)) < 1
\ee
and $\kappa_n (q_1 - q_0) \geq 0$ near infinity.
\end{theorem}

\begin{proof}
Lemma~\ref{lem:HS} together with (\ref{estgfper}) shows, that Hypothesis~\ref{hyp:h0h1}
is satisfied if $q_0 - q_1 \in L^1(1,\infty)$. Thus we can
either apply Theorem~\ref{thm:wronskzeros} or Theorem~\ref{thmsing}
to conclude that $\tau_1-\lam$ is nonoscillatory with respect to
$\tau_0-\lam$ for any $\lam\in\R\backslash\sig_{ess}(H_0)$.
Hence, Theorem~\ref{thm:rosc} (ii) is applicable and it suffices to show that
$\tau_1$ is relatively oscillatory (resp.\ nonoscillatory) with respect to $\tau_0$.

Without restriction, we assume $\kappa_n > 0$.
For the first statement, note that we can find $c$, $\eps>0$ such that 
$$
q_0(x) < q_0(x) + \frac{1}{(\kappa_n - \eps)x^2} < q_1(x), \qquad x>c.
$$
Now since, perturbations with compact support only
add finitely many eigenvalues, it is no restriction to assume $c = 1$.

Next, \cite[Thm.~1]{kms} shows that $\tau_0 + (\kappa_n-\eps)^{-1} x^{-2} - E_n$
is relatively oscillatory with respect to $\tau_0 - E_n$. Hence $\tau_1 - E_n$ being
relatively oscillatory with respect to $\tau_0 - E_n$ now follows using
our comparison theorem for Wronskians (cf.\ Theorem~\ref{thm:scw}).

For the second statement, we first note that our conditions imply
$$
q_0(x) \leq q_1(x) < q_0(x) + \frac{1}{(\kappa_n+\eps)x^2}
$$
near infinity, and then one proceeds as before.
\end{proof}

\noindent
We note that even the second order term was also computed in \cite{kms}. So
we obtain:

\begin{theorem}
Assume 
\be
\lim_{x\rightarrow \infty} \kappa_n x^2 (q_1 (x) - q_0(x)) = 1
\ee
in addition to the assumptions in Theorem~\ref{thm:app}.
Then $E_n$ is an accumulation point of eigenvalues if
\be
\liminf_{x \to\infty} \log^2(x) (\kappa_n x^2 (q_1 (x) - q_0(x)) -1) > 1,
\ee
and $E_n$ is not an accumulation point of eigenvalues if
\be
\limsup_{x \to\infty} \log^2(x) (\kappa_n x^2 (q_1 (x) - q_0(x)) -1) < 1.
\ee
\end{theorem}

\begin{proof}
Similarly to that of Theorem~\ref{thm:app}, except, one now
uses \cite[Thm.~2]{kms}.
\end{proof}

\noindent
The main argument in \cite{kms} is a perturbation argument for the difference
of Pr\"ufer angles for the solution of the unperturbed and the
perturbed equation. This corresponds exactly to calculating the
asymptotics of the Pr\"ufer angle of the Wronskian.

\section{More on Pr\"ufer angles and the case of regular operators}
\label{sec:proofreg}

Now let us suppose that $\tau_{0,1}$ are both regular at $a$ and $b$ with boundary conditions
\be \label{bc}
\cos(\alpha) f(a) - \sin(\alpha) p(a)f'(a) =0 , \quad
\cos(\beta) f(b) - \sin(\beta) p(b)f'(b) =0.
\ee
Hence we can choose $\psi_\pm(\lam,x)$ such that
$\psi_-(\lam,a)=\sin(\alpha)$, $p(a)\psi_-'(\lam,a)=\cos(\alpha)$ respectively
$\psi_+(\lam,b)=\sin(\beta)$, $p(b)\psi_+'(\lam,b)=\cos(\beta)$. In particular,
we may choose
\be \label{normtha}
\theta_-(\lam,a) =\alpha \in [0,\pi), \quad -\theta_+(\lam,b) = \pi -\beta \in
[0,\pi).
\ee

Next we introduce
\be
\tau_\eps = \tau_0 + \frac{\eps}{r} (q_1 -q_0)
\ee
and investigate the dependence with respect to $\eps\in[0,1]$. If $u_\eps$
solves $\tau_\eps u_\eps = 0$, then the corresponding Pr\"ufer angles satisfy
\be
\dot{\theta}_\eps(x) = -\frac{W_x(u_\eps, \dot{u}_\eps)}{\rho_\eps^2(x)},
\ee
where the dot denotes a derivative with respect to $\eps$.

\begin{lemma} \label{prwpsiepsdot}
We have
\be
W_x(\psi_{\eps,\pm}, \dot{\psi}_{\eps,\pm}) = \left\{
\begin{array}{l} \int_x^b (q_0(t)-q_1(t)) \psi_{\eps,+}(t)^2 dt \\ -\int_a^x
(q_0(t)-q_1(t)) \psi_{\eps,-}(t)^2 dt \end{array} \right. ,
\ee
where the dot denotes a derivative with respect to $\eps$ and
$\psi_{\eps,\pm}(x)= \psi_{\eps,\pm}(0,x)$.
\end{lemma}

\begin{proof}
Integrating (\ref{dwr}) we obtain
\be
W_x(\psi_{\eps,\pm}, \psi_{\ti{\eps},\pm}) = (\ti{\eps}-\eps) \left\{
\begin{array}{l} \int_x^b (q_0(t)-q_1(t)) \psi_{\eps,+}(t) \psi_{\ti{\eps},+}(t) dt, \\
- \int_a^x (q_0(t)-q_1(t)) \psi_{\eps,-}(t) \psi_{\ti{\eps},-}(,t) dt. \end{array} \right.
\ee
Now use this to evaluate the limit
\be
\lim_{\ti{\eps} \to \eps}W_x \Big(\psi_{\eps,\pm},
\frac{\psi_{\pm,\eps} - \psi_{\ti{\eps},\pm}}{\eps-\ti{\eps}} \Big).
\ee
\end{proof}

\noindent
Denoting the Pr\"ufer angles of $\psi_{\eps,\pm}(x)= \psi_{\eps,\pm}(0,x)$ by $\theta_{\eps,+}(x)$,
this result implies for $q_0-q_1\ge 0$,
\begin{align} \nn
\dot{\theta}_{\eps,+}(x) &= -\frac{\int_x^b (q_0(t)-q_1(t)) \psi_{\eps,+}(t)^2 dt}{\rho_{\eps,+}(x)^2}
\le 0, \\ \label{thetadot}
\dot{\theta}_{\eps,-}(x) &= \frac{\int_a^x (q_0(t)-q_1(t)) \psi_{\eps,-}(t)^2
dt}{\rho_{\eps,-}(x)^2} \ge 0,
\end{align}
with strict inequalities if $q_1 \not\equiv q_0$.

Now we are ready to investigate the associated operators $H_0$ and $H_1$.
In addition, we will choose the same boundary conditions for $H_\eps$ as
for $H_0$ and $H_1$.

\begin{lemma} \label{lemevheps}
Suppose $q_0-q_1\ge 0$ (resp $q_0-q_1\leq 0$).
Then the eigenvalues of $H_\eps$ are analytic functions with respect to
$\eps$ and they are decreasing (resp.\ increasing).
\end{lemma}

\begin{proof}
First of all the Pr\"ufer angles $\theta_{\eps,\pm}(x)$
are analytic with respect to $\eps$ since
$\tau_\eps$ is by a well-known result from
ordinary differential equations (see e.g., \cite[Thm.~13.III]{wa}). Moreover,
$\lam\in\sig(H_\eps)$ is equivalent to $\theta_{\eps,+}(a)\equiv \alpha \mod \pi$
(respectively $\theta_{\eps,-}(b) \equiv \beta \mod \pi$),
where $\alpha$ (respectively $\beta$) generates the
boundary condition (cf.\ (\ref{bc})).
\end{proof}

\noindent
In particular, this implies that $\dim\Ran P_{(-\infty,\lam)} (H_\eps)$ is continuous
from below (resp. above) in $\eps$ if $q_0-q_1\ge 0$ (resp $q_0-q_1\leq 0$).

Now we are ready for the

\begin{proof}[Proof of Theorem~\ref{thm:reg}]
It suffices to prove the result for $\#(\psi_{0,+}, \psi_{\eps,-})$.
Again we can assume $\lam_0=\lam_1=0$ without restriction. 
We split $q_0-q_1$ according to
$$
q_0 - q_1 = q_+ - q_-,\qquad q_+, q_- \geq 0,
$$
and introduce the operator $\tau_- = \tau_0 - q_-/r$.
Then $\tau_-$ is a negative perturbation of $\tau_0$ and
$\tau_1$ is a positive perturbation of $\tau_-$.

Furthermore define $\tau_\eps$ by
$$
\tau_\eps = \begin{cases}
\tau_0 + 2\eps (\tau_- - \tau_0), & \eps\in[0,1/2],\\
\tau_- + 2(\eps -1/2)(\tau_1 - \tau_-), & \eps\in[1/2,1].
\end{cases}
$$
Let us look at
$$
N(\eps)=\#(\psi_{0,+}, \psi_{\eps,-}) =
\ceil{\Delta_\eps(b)/\pi} - \floor{\Delta_\eps(a)/\pi} -1, \quad
\Delta_\eps(x)=\Delta_{\psi_{0,+}, \psi_{\eps,- }}(x)
$$
and consider $\eps\in [0,1/2]$.
At the left boundary $\Delta_\eps(a)$ remains constant whereas at the right
boundary $\Delta_\eps(b)$ is increasing by Lemma~\ref{prwpsiepsdot}.
Moreover, it hits a multiple of $\pi$ whenever $0\in\sig(H_\eps)$.
So $N(\eps)$ is a piecewise constant function which is continuous from below
and jumps by one whenever $0\in\sig(H_\eps)$. By Lemma~\ref{lemevheps}
the same is true for
$$
P(\eps) = \dim\Ran\, P_{(-\infty, 0)} (H_\eps) - \dim\Ran\,
P_{(-\infty, 0]} (H_0)
$$
and since we have $N(0)=P(0)$, we conclude $N(\eps)=P(\eps)$ for all
$\eps\in[0,1/2]$. To see the remaining case $\eps=[1/2,1]$, simply
replace increasing by decreasing and continuous from below by continuous
from above.
\end{proof}

\section{Approximation by regular operators}
\label{sec:appreg}

Now we want to extend our results to singular operators. We will do so by
approximating a singular operator by a sequence of regular ones following
\cite[Chap.~14]{wdln}.

Abbreviate in the following $L^2((c,d),r\,dx)$ as $L^2(c,d)$.
Fix functions $u,v\in\dom(\tau)$ and pick $c_n\downarrow a$,
$d_n\uparrow b$. Define $\ti{H}_n$
\be \label{deftiHm}
\ti{H}_n: \begin{array}[t]{lcl} \dom(\ti{H}_n) &\to& L^2(c_n,d_n) \\ f
&\mapsto& \tau f \end{array},
\ee
where
\be
\dom(\tilde{H}_n) = \{ f \in L^2(c_n,d_n) | \begin{array}[t]{l} f, pf' \in
AC(c_n,d_n), \, \tau f \in L^2(c_n,d_n),\\ W_{c_n}(u,f) =
W_{d_n}(v,f) =0 \}.
\end{array}
\ee
Take $H_n=\alpha \id \oplus \ti{H}_n \oplus\alpha \id$ on
$L^2(a,b)=L^2(a,c_n)\oplus L^2(c_n,d_n)\oplus L^2(d_n, b)$, where
$\alpha$ is a fixed real constant. Then we have the following result:

\begin{lemma} \label{strconh}
Suppose that either $H$ is limit point at $a$ or
that $u = \psi_-(\lam_0)$ for some $\lam_0\in\R$ and similarly, that either
$H$ is limit point at $b$ or $v=\psi_+(\lam_1)$ for some $\lam_1\in\R$.
Then $H_n$ converges to $H$ in strong resolvent sense as $n\to\infty$.

Furthermore, if $H$ is limit circle at $a$ (resp. $b$), we can
replace $u$ (resp. $v$) with any function in $\dom(\tau)$, which
generates the boundary condition.
\end{lemma}

\noindent
However, strong resolvent convergence is not sufficient for our purpose here. In addition
we will need the following result from \cite{sw} (see also \cite{wdln}).
We give a slightly refined analysis which allows eigenvalues at the boundary of
the spectral intervals and the possibility of infinite-dimensional projections.
We remark that for a self-adjoint projector $P$ we have
\be
\dim\Ran(P) = \tr(P) = \|P\|_{\mathcal{J}^1},
\ee
where $\|.\|_{\mathcal{J}^1}$ denotes the trace class norm. If $P$ is not finite-rank,
all three numbers equal $\infty$.
Then we have the following result (\cite[Lem.~2]{te}, see also \cite{sw}):

\begin{lemma}\label{strecon}
Let $A_n\to A$ in strong resolvent sense and suppose
$\tr(P_{(\lam_0,\lam_1)}(A_n)) \le \tr(P_{(\lam_0,\lam_1)}(A))$.

Then,
\be
\lim_{n\to\infty} \tr(P_{(\lam_0,\lam_1)}(A_n)) = \tr(P_{(\lam_0,\lam_1)}(A)),
\ee
and if $\tr(P_{(\lam_0,\lam_1)}(A))<\infty$, we have
\be
\lim_{n\to\infty} \| P_{(\lam_0,\lam_1)}(A_n) - P_{(\lam_0,\lam_1)}(A)\|_{\mathcal{J}^1} =0.
\ee
\end{lemma}

\begin{proof}
This follows from (see e.g.\ \cite[Lem.~5.2]{gst})
\be\label{liminfsrc}
\tr(P_{(\lam_0,\lam_1)}(A)) \le \liminf_{n\to\infty} \tr(P_{(\lam_0,\lam_1)}(A_n)),
\ee
together with Gr\"umm's theorem (\cite[Thm.~2.19]{str}).
\end{proof}

\noindent
The key result of Stolz and Weidmann is that this lemma is applicable
if certain Weyl solutions $\psi_\pm(\lam)$ are used to generate the
boundary conditions of $\ti{H}_n$. As already pointed out, the
version below is slightly refined since it allows $\lam_0$, $\lam_1$ to be
eigenvalues of $H$.

\begin{lemma}[\cite{sw}]\label{lem:swres0}
Suppose $[\lam_0, \lam_1]\cap\sigma_{ess}(H) = \emptyset$ and let
$H_n$ be defined as in (\ref{deftiHm}) with $u = \psi_-(\lam_-)$,
$v=\psi_+(\lam_+)$ and $\lam_\pm\in[\lam_0,\lam_1]$. Then, 
\be
\tr(P_{(\lam_0, \lam_1)}(\ti{H}_n)) \le  \tr(P_{(\lam_0, \lam_1)}(H)).
\ee
\end{lemma}

\begin{proof}
Abbreviate $P=P_{(\lam_0, \lam_1)}(H)$, $P_n=P_{(\lam_0, \lam_1)}(\ti{H}_n)$.
For $\ti{\psi}_1, \dots, \ti{\psi}_k \in \Ran P_n$ being eigenfunctions of $\ti{H}_n$, construct
$$
\psi_j (x) = \left\{\begin{array}{cl} \gam_{j,u} u(x), & x < c_n,\\
\ti{\psi}_j(x), & c_n \leq x \leq d_n,\\
\gam_{j,v} v(x), & x > d_n, \end{array} \right.
$$
where $\gam_{j,u}$, $\gam_{j,v}$ are chosen such that $\psi_j$ and $p \psi_j'$ are continuous.
A computation shows that
$$
\|(H - \frac{\lam_1 + \lam_0}{2}) \psi\| < \frac{\lam_1 - \lam_0}{2} \|\psi\|
$$
for any $\psi$ in the linear span of the $\psi_j$'s, which yields the first result.
\end{proof}

\noindent
This version is sufficient to give an alternative proof for the main theorem in \cite{gst}.

\begin{proof}[Proof of Theorem~\ref{thm:gst}]
Approximate $H_0$ by regular operators $H_n$ defined as in (\ref{deftiHm}) with
$u=\psi_{0,-}(\lam_0)$, $v=\psi_{0,+}(\lam_1)$. Denote by $\psi_{0,\pm}^n(\lam)$ the solutions
of the approximating problems. Then, by construction $\psi_{0,-}^n(\lam_0,x)=\psi_{0,-}(\lam_0,x)$
respectively $\psi_{0,+}^n(\lam_1,x)=\psi_{0,+}(\lam_1,x)$ for $x\in(c_n,d_n)$ and
Theorem~\ref{thm:reg} in the special case $H_1=H_0$ implies
$$
 \tr(P_{(\lam_0, \lam_1)}(\ti{H}_n)) = \#_{(c_n,d_n)}(\psi_{0,-}(\lam_0),\psi_{0,+}(\lam_1)).
$$
Letting $n\to\infty$ the left-hand side converges to $\tr(P_{(\lam_0, \lam_1)}(H_0))$ by the first
part of Lemma~\ref{lem:swres0}. Hence the right-hand side converges as well and,
according to Definition~\ref{def:wsf}, is given by $\#(\psi_{0,\mp}(\lam_0),\psi_{0,\pm}(\lam_1))$.
\end{proof}

\noindent
However, the proof of our Theorem~\ref{thm:wronskzeros} requires some further extensions.
In fact, in \cite{sw} Stolz and Weidmann point out that the Weyl functions of a different operator
$\ti{H}$ will also do, as long as $\ti{H}$ is not too far away from $H$. Again
the version below is slightly improved to allow for some border line cases.

\begin{lemma}[\cite{sw}]\label{lem:swres1}
Suppose $[\lam_0,\lam_1]\cap\sig_{ess}(H)=\emptyset$.
Let $\ti{\tau} = \tau + \ti{q}$, where $\ti{q}$ is bounded,
and pick the same boundary conditions for $\ti{H}$ as for $H$ (if any).
Abbreviate
\be
Q_a = [\liminf_{x\to  a} \ti{q}(x), \limsup_{x\to  a} \ti{q}(x)], \quad
Q_b = [\liminf_{x\to  b}\ti{q}(x), \limsup_{x\to  b} \ti{q}(x)],
\ee
and choose $\lam_-$ such that one of following conditions holds:
\begin{enumerate}
\item $\lam_- - Q_a \subseteq (\lam_0,\lam_1)$, or
\item $\lam_- - Q_a \subseteq [\lam_0,\lam_1)$ and $\ti{q} \le 0$ near $a$, or
\item $\lam_- - Q_a \subseteq (\lam_0,\lam_1]$ and $\ti{q} \ge 0$ near $a$.
\end{enumerate}
Similarly, choose $\lam_+$ to satisfy one of these conditions with $a$ replaced by $b$.

Then, $H_n$ defined as in (\ref{deftiHm}) with $u = \ti{\psi}_-(\lam_-)$, $v=\ti{\psi}_+(\lam_+)$,
satisfies
\be
\limsup_{n\to\infty} \tr(P_{(\lam_0, \lam_1)}(\ti{H}_n)) \le  \tr(P_{(\lam_0, \lam_1)}(H)).
\ee

Furthermore, if we just require
\be
\lam_- - Q_a \subseteq [\lam_0,\lam_1], \qquad \lam_+ - Q_b \subseteq [\lam_0,\lam_1],
\ee
we at least have
\be
\limsup_{n\to\infty} \tr(P_{[\lam_0, \lam_1]}(\ti{H}_n)) \le  \tr(P_{[\lam_0, \lam_1]}(H)).
\ee
\end{lemma}

\begin{proof}
Since any of our three conditions implies
$|\lam_- -\frac{\lam_1+\lam_0}{2} - \ti{q}(x)| \le \frac{\lam_1-\lam_0}{2}$ for $x$ sufficiently
close to $a$ and similarly $|\lam_+ -\frac{\lam_1+\lam_0}{2} - \ti{q}(x)| \le \frac{\lam_1-\lam_0}{2}$
for $x$ sufficiently close to $b$, we can proceed as in the previous lemma to prove the
first claim.

For the second claim choose $n$ sufficiently large such that
$|\lam_\pm - \ti{q}(x) - 2^{-1}(\lam_1+\lam_0)| \le (2^{-1}(\lam_1-\lam_0) + \eps)$ for $x<c_n$,
respectively, $x>d_n$. Then, with the same argument as in the previous lemma, we have
$$
\|(H - \frac{\lam_1 + \lam_0}{2}) \psi\| \le (\frac{\lam_1 - \lam_0}{2}+ \eps) \|\psi\|
$$
and hence the second claim follows.
\end{proof}

\noindent
Since our results involve projections to half-open intervals, we need one
further step.

\begin{lemma}\label{lem:swres2}
Suppose $[\lam_0,\lam_1]\cap\sig_{ess}(H)=\emptyset$.
Let $\ti{\tau} = \tau + \ti{q}$, where $\lim_{x\to a} \ti{q}(x)=0$ if $a$ is singular and
$\lim_{x\to b} \ti{q}(x)=0$ if $b$ is singular.
Furthermore, pick the same boundary conditions for $\ti{H}$ as for $H$ (if any).

Define $H_n$ as in (\ref{deftiHm}) with $u = \ti{\psi}_-(\lam)$, $v=\ti{\psi}_+(\lam)$,
$\lam\in\{\lam_1,\lam_2\}$. If $\lam=\lam_1$ and $\ti{q}\leq 0$ (near $a$ and $b$), then
\be
\lim_{n\to\infty} \tr(P_{(\lam_0,\lam_1]}(\ti{H}_n)) = \tr(P_{(\lam_0,\lam_1]}(H)),
\ee
and if $\lam=\lam_0$ and $\ti{q}\geq 0$, then
\be
\lim_{n\to\infty} \tr(P_{[\lam_0,\lam_1)}(\ti{H}_n)) = \tr(P_{[\lam_0,\lam_1)}(H)).
\ee
\end{lemma}

\begin{proof}
Without restriction, we just proof the first claim.
For a sufficiently small $\eps>0$ we still have
$[\lam_0,\lam_1+\eps]\cap\sig_{ess}(H)=\emptyset$ and
thus by the previous lemma
$$
\lim_{n\to\infty} \tr(P_{(\lam_0,\lam_1+\eps)}(\ti{H}_n)) =
\tr(P_{(\lam_0,\lam_1+\eps)}(H))
$$
and
$$
\lim_{n\to\infty} \tr(P_{(\lam_1,\lam_1+\eps)}(\ti{H}_n)) =
\tr(P_{(\lam_1,\lam_1+\eps)}(H)).
$$
Hence the result follows from
$P_{(\lam_0,\lam_1]} = P_{(\lam_0,\lam_1+\eps)}-P_{(\lam_1,\lam_1+\eps)}$.
\end{proof}

\noindent
Finally, we note:

\begin{lemma} \label{lem:csf}
Suppose $\tau_j u_j = \lam_j u_j$, $j=0,1$, with $q_1\le q_0$ near singular endpoints and
$\lam_0<\lam_1$. If $\tau_j u_{j,n} = \lam_{j,n} u_{j,n}$, where
$\lam_{j,n} \to \lam_j$ and $u_{j,n} \to u_j$, uniformly on compact sets $[c,d]\subseteq(a,b)$, then
\be
\liminf_{n\to\infty} \#(u_{0,n},u_{1,n}) \geq \#(u_0,u_1).
\ee
\end{lemma}

\begin{proof}
Let $N\in\N_0$ be any finite number with $N \le \#(u_0,u_1)$. Choose a compact set $[c,d]$
containing $N$ sign flips of $W(u_0,u_1)$. Then, for $n$ sufficiently large,
$W(u_{0,n},u_{1,n})$ has $N$ sign flips in $[c,d]$. Hence
$\#(u_{0,n},u_{1,n}) \geq \#_{(c,d)}(u_{0,n},u_{1,n}) = N$ and the claim follows.
\end{proof}

\noindent
Now, we are ready for the

\begin{proof}[Proof of Theorem~\ref{thm:wronskzeros}]
It suffices to show the $\#(\psi_{1,-}(\lam_j), \psi_{0,+}(\lam_j))$ case.
Define $\ti{H}_{j,n}$, $j=0,1$, as in (\ref{deftiHm}) with $u=\psi_{1,-}(\lam_0)$ and
$v=\psi_{0,+}(\lam_0)$.

Denote by $\psi_{j,\pm}^n(\lam)$, $j=0,1$, the solutions of the approximating problems.
Then, by Theorem~\ref{thm:reg},
$$
\tr(P_{(-\infty,\lam_0)}(\ti{H}_{1,n})) - \tr(P_{(-\infty,\lam_0]}(\ti{H}_{0,n})) =
\#_{(c_n,d_n)}(\psi_{1,-}^n(\lam_0),\psi_{0,+}^n(\lam_0))
$$
and we need to investigate the limits as $n\to\infty$.

First of all $\psi_{1,-}^n(\lam_0,x)= \psi_{1,-}(\lam_0,x)$, $\psi_{0,+}^n(\lam_0,x)= \psi_{0,+}(\lam_0,x)$
for $x\in(c_n,d_n)$ implies
\begin{align*}
\lim_{n\to\infty} \#_{(c_n,d_n)}(\psi_{1,-}^n(\lam_0),\psi_{0,+}^n(\lam_0))
&= \lim_{n\to\infty} \#_{(c_n,d_n)}(\psi_{1,-}(\lam_0),\psi_{0,+}(\lam_0))\\
&= \#(\psi_{1,-}(\lam_0),\psi_{0,+}(\lam_0)).
\end{align*}
This takes care of the number of sign flips and it remains to look at the spectral
projections. Let $\lam_0<\sig_{ess}(H_0)$, that is $H_0$ and hence also $H_1$
are bounded from below. Replacing $P_{(-\infty,\lam_0)}(H_j)$ by $P_{(\lam,\lam_0)}(H_j)$
with some $\lam$ below the spectrum of both $H_0$ and $H_1$ we infer from
Lemma~\ref{lem:swres2}
$$
\lim_{n\to\infty} \tr(P_{(-\infty,\lam_0)}(\ti{H}_{1,n})) = \tr(P_{(-\infty,\lam_0)}(H_1))
$$
and
$$
\lim_{n\to\infty} \tr(P_{(-\infty,\lam_0]}(\ti{H}_{0,n})) = \tr(P_{(-\infty,\lam_0]}(H_0)).
$$
This settles the first claim (\ref{eq:wronskzeros0}), where $\lam_0<\sig_{ess}(H_0)$.

For the second claim (\ref{eq:wronskzeros1}), we first note that $\tau_1-\lam_0$ is relatively
nonoscillatory with respect to $\tau_0-\lam_0$ by Lemma~\ref{lem:nonoscingap}.
Next note that $\psi_{0,+}^n(\lam_1,.) \rightarrow \psi_{0,+}(\lam_1,.)$ pointwise, since
\be
\psi_{0,+}^n(\lam_1,x) = c_0(\lam_1,x) + m_{0,+}^n(\lam_1) s_0(\lam_1,x),
\ee
where $c_0(\lam,x)$, $s_0(\lam,x)$ is a fundamental system of solutions for $\tau_0-\lam$, and
$m_{0,+}^n(\lam)$ are the corresponding Weyl--Titchmarsh $m$-functions. Next,
strong resolvent convergence implies convergence of the Weyl
$m$-function and hence uniform convergence of
$\psi_{0,+}^n(\lam_1,x)\to \psi_{0,+}(\lam_1,x)$ on compact sets. Clearly the same applies
to $\psi_{1,-}^n(\lam_1,x) \to \psi_{1,-}(\lam_1,x)$.
Thus, by Lemma~\ref{strecon}, Lemma~\ref{lem:swres2}, and Lemma~\ref{lem:csf},
\begin{align}
\label{eq:wrprin1}
\tr(P_{[\lam_0, \lam_1)}(H_1)) &- \tr(P_{(\lam_0, \lam_1]}(H_0))\\ \nn
&\geq \#(\psi_{1,-} (\lam_1), \psi_{0,+} (\lam_1)) - 
\#(\psi_{1,-} (\lam_0), \psi_{0,+} (\lam_0)).
\end{align}
Repeating the argument with $u=\psi_{1,-}(\lambda_1)$ and $v=\psi_{0,+}(\lambda_1)$
shows that
\begin{align}
\label{eq:wrprin2}
\tr(P_{[\lam_0, \lam_1)}(H_1)) &- \tr(P_{(\lam_0, \lam_1]}(H_0))\\ \nn
&\leq \#(\psi_{1,-} (\lam_1), \psi_{0,+} (\lam_1)) - 
\#(\psi_{1,-} (\lam_0), \psi_{0,+} (\lam_0)).
\end{align}
This proves the second claim.
\end{proof}

\section{Approximation in trace norm}
\label{sec:apptr}

Now we begin with an alternative approach toward singular differential
operators by proving the case where
$q_1-q_0$ has compact support. The next lemma would follow from
Theorem~\ref{thm:wronskzeros}, but to demonstrate that this
approach is independent of the last, we will provide an alternative
proof.

\begin{lemma} \label{lem:compactcase}
Let $H_j$, $j=0,1$, be Sturm--Liouville operators on $(a,b)$ associated with $\tau_j$, and
suppose that $r^{-1}(q_1-q_0)$ has support in a bounded interval $(c,d)\subseteq (a,b)$, where
$a<c$ if $a$ is singular and $d<b$ if $b$ is singular. Moreover, suppose $H_0$ and $H_1$
have the same boundary conditions (if any).

Suppose $\lam_0 < \inf\sig_{ess}(H_0)$. Then,
\be
\dim\Ran P_{(-\infty,\lam_0)}(H_1) - \dim\Ran P_{(-\infty,\lam_0]}(H_0)
= \#(\psi_{1,\mp} (\lam_0), \psi_{0,\pm} (\lam_0)).
\ee

Suppose $\sig_{ess}(H_0) \cap [\lam_0,\lam_1] = \emptyset$. Then,
\begin{align}
\nonumber
\dim\Ran P_{[\lam_0, \lam_1)} (H_1) &- \dim\Ran P_{(\lam_0, \lam_1]} (H_0)\\
&= \#(\psi_{1,\mp} (\lam_1), \psi_{0,\pm} (\lam_1)) - 
\#(\psi_{1,\mp} (\lam_0), \psi_{0,\pm} (\lam_0)).
\end{align}
\end{lemma}

\begin{proof}
By splitting $r^{-1}(q_1-q_0)$ into a positive and negative part as in the
proof of the regular case (Theorem~\ref{thm:reg}), we can reduce it to the
case where $r^{-1}(q_1-q_0)$ is of one sign, say $r^{-1}(q_1-q_0)\ge 0$.
Define $H_\eps= \eps H_1 + (1-\eps) H_0$ and observe that
$\psi_{\eps,-}(z,x)=\psi_{0,-}(z,x)$ for $x\le c$, respectively,
$\psi_{\eps,+}(z,x)=\psi_{0,+}(z,x)$ for $x\ge d$. Furthermore,
$\psi_{\eps,\pm}(z,x)$ is analytic with respect to $\eps$ and
$\lam\in\sig_p(H_\eps)$ if and only if $W_d(\psi_{0,+}(\lam),\psi_{\eps,-}(\lam))=0$.
Now the proof can be done as in the regular case. 
\end{proof}

\begin{lemma} \label{lem:compactxi}
Suppose $H_0$, $H_1$ satisfy the same assumptions as in the previous lemma and
set $H_\eps= \eps H_1 + (1-\eps) H_0$. Then,
\be
\| \sqrt{r^{-1}|q_0-q_1|} R_{H_\eps}(z) \|_{\mathcal{J}_2} \le C(z), \qquad \eps\in[0,1].
\ee
In particular, the resolvent difference of $H_0$ and $H_1$ is trace class and
\be
\xi(\lambda, H_1, H_0) =  \#(\psi_{1,\mp}(\lam), \psi_{0,\pm}(\lam))
\ee
for every $\lam \in\R\cap\rho(H_0)\cap\rho(H_1)$.
Here $\xi(H_1,H_0)$ is assumed to be constructed
such that $\eps \mapsto \xi(H_\eps, H_0)$ is a continuous mapping from $[0,1] \rightarrow
L^1(\R,(\lam^2+1)^{-1} d\lam)$.
\end{lemma}

\begin{proof}
Denote by
$$
G_\eps(z,x,y) = (H_\eps-z)^{-1}(x,y) = \frac{\psi_{\eps,-}(z,x_<),\psi_{\eps,+}(z,y_>)}{W(\psi_{\eps,-}(z),\psi_{\eps,+}(z))},
$$
where $x_< = \min(x,y)$, $y_>=\max(x,y)$, the Green's function of $H_\eps$. As pointed out
in the proof of the previous lemma, $\psi_{\eps,\pm}(z,x)$ is analytic with respect to $\eps$
and hence a simple estimate shows
$$
\int_a^b \int_a^b |G_\eps(z,x,y)|^2 |r(y)^{-1}(q_1(y)-q_0(y))| r(x)dx\, r(y)dy \le C(z)^2
$$
for $\eps\in[0,1]$, which establishes the first claim.

Furthermore, a straightforward calculation (using (\ref{dwr})) shows
$$
G_{\eps'}(z,x,y) = G_\eps(z,x,y) +
(\eps-\eps') \int_a^b G_{\eps'}(z,x,t) r^{-1}(t)(q_1(t)-q_0(t)) G_\eps(z,t,y) r(t) dt.
$$
(Note that this does not follow from the second resolvent identity unless $r^{-1}(q_1-q_0)$ is
relatively bounded with respect to $H_0$.) Hence, $R_{H_{\eps'}}(z) - R_{H_\eps}(z)$
can be written as the product of two Hilbert--Schmidt operators whose norm can
be estimated by the first claim,
\be
\| R_{H_{\eps'}}(z) - R_{H_\eps}(z) \|_{\mathcal{J}_1} \le |\eps'-\eps| C(z)^2.
\ee
Thus $\eps \mapsto \xi(H_\eps, H_0)$ is continuous by Lemma~\ref{lem:Heps}.
The rest follows from (\ref{dimxi}).
\end{proof}

\noindent
Now we come to the

\begin{proof}[Proof of Theorem~\ref{thmsing}]
We first assume that we have compact support near one endpoint, say $a$.
Furthermore, abbreviate $V=r^{-1}(q_0-q_1)$ which satisfies Hypothesis~\ref{hyp:h0h1} 
(and thus also Hypothesis~\ref{hyp:h0v}).

Define by $K_\eps$ the operator of multiplication by $\chi_{(a,b_\eps]}$ with $b_\eps \uparrow b$
as $\eps\uparrow 1$. Then $K_\eps$ satisfies the assumptions of Lemma~\ref{lem:resolvconv}.
Introduce $H_\eps = H_0 - K_\eps V$, and denote by $\psi_{\eps,-}(\lam,x)$ the corresponding
solutions satisfying the boundary condition at $a$.

By Lemma~\ref{lem:resolvconv} we have $\xi(., H_\eps, H_0) \to \xi(., H_1, H_0)$ as
$\eps\to 1$ in $L^1(\R,(\lam^2+1)^{-1}d\lam)$. Moreover,  $H_\eps \to H_1$ in (trace) norm
resolvent sense and hence $\lam\in\rho(H_1)$ implies $\lam\in\rho(H_\eps)$ for $\eps$
sufficiently close to $1$. Since $\xi(\lam, H_\eps, H_0) \in \Z$ is constant near every
$\lam\in\R\cap\rho(H_0)\cap\rho(H_\eps)$, we must have
$\xi(\lam, H_\eps, H_0)= \xi(\lam, H_1, H_0)$ for $\eps \ge \eps_0$ with some
$\eps_0$ sufficiently close to $1$.

Now let us turn to the Wronskians. We first prove the $\#(\psi_{1,-}(\lam), \psi_{0,+}(\lam))$ case.
By Lemma~\ref{lem:compactxi} we know
$$
\xi(\lambda, H_\eps, H_0) = \#(\psi_{\eps,-}(\lam),\psi_{0,+}(\lam)
$$
for every $\eps<1$. Concerning the right-hand side observe that
$$
W_x(\psi_{\eps,-}(\lam), \psi_{0,+}(\lam)) = W_x(\psi_{1,-}(\lam), \psi_{0,+}(\lam))
$$
for $x \le b_\eps$ and that $W_x(\psi_{\eps,-}(\lam), \psi_{0,+}(\lam))$ is constant for
$x \ge b_\eps$. This implies that for $\eps\ge \eps_0$ we have
\begin{align*}
\xi(\lam, H_1, H_0) &= \xi(\lam, H_\eps, H_0) = \#(\psi_{\eps,-}(\lam), \psi_{0,+}(\lam))\\
& = \#_{(a,b_\eps)}(\psi_{\eps,-}(\lam), \psi_{0,+}(\lam)) = \#_{(a,b_\eps)}(\psi_{1,-}(\lam), \psi_{0,+}(\lam)).
\end{align*}
In particular, the last item $\#_{(a,b_\eps)}(\psi_{1,-}(\lam), \psi_{0,+}(\lam))$ is eventually constant
and thus has a limit which, by Definition~\ref{def:wsf}, is $\#(\psi_{1,-}(\lam), \psi_{0,+}(\lam))$.

For the corresponding $\#(\psi_{1,+}(\lam), \psi_{0,-}(\lam))$ case one simply exchanges the
roles of $H_0$ and $H_1$.

Hence the result holds if the perturbation has compact support near one endpoint. Now one repeats the
argument to remove the compact support assumption near the other endpoint as well.
\end{proof}

\section{Appendix: Some facts on the spectral shift function}
\label{sec:xi}

In this appendix we collect some facts on Krein's spectral shift function which are of
relevance to us. Most results are taken from \cite{yafams} (see also \cite{sro} for
an easy introduction).

Two operators $H_0$ and $H_1$ are called resolvent comparable, if
\be\label{eq:resolvcond}
R_{H_1}(z) - R_{H_0}(z)
\ee
is trace class for one $z\in\rho(H_1)\cap\rho(H_0)$. By the first resolvent identity \eqref{eq:resolvcond}
then holds for all $z\in\rho(H_1)\cap\rho(H_0)$.

\begin{theorem}[Krein \cite{krein}]\label{thm:shifttheorem}
Let $H_1$ and $H_0$ be two resolvent comparable self-adjoint operators, 
then there exists a function
\be
\xi(\lam, H_1, H_0) \in L^1(\R, (\lambda^2 + 1)^{-1}d\lam)
\ee
such that 
\be\label{eq:traceformula}
\tr(f(H_1) - f(H_0)) = \int_{-\infty}^\infty \xi(\lam, H_1, H_0) f'(\lambda) d\lam
\ee
for every smooth function $f$ with compact support.
\end{theorem}
Note: Equation \eqref{eq:traceformula} holds in fact for a much larger class of functions $f$.
See \cite[Thm.~8.7.1]{yafams} for this and a proof of the last theorem.

The function $\xi(\lam) = \xi(\lam, H_1, H_0)$ is called Krein's spectral shift function and is
unique up to an additive constant. Moreover, $\xi(\lam)$ is constant on every interval
$(\lam_0, \lam_1) \subset \rho(H_0)\cap\rho(H_1)$. Hence, if
$\dim\Ran P_{(\lam_0, \lam_1)}(H_j) < \infty$, $j = 0,1$,
then $\xi(\lam)$ is a step function and
\be \label{dimxi}
\dim \Ran P_{(\lam_0, \lam_1)} (H_1) - \dim \Ran P_{(\lam_0, \lam_1)} (H_0) =
\lim_{\eps\downarrow0} \Big( \xi(\lam_1-\eps) - \xi(\lam_0+\eps) \Big).
\ee
This formula explains the name spectral shift function.

Before investigating further the properties of the SSF,
we will recall a few things about trace ideals
(see for example \cite{str}). First,
for $1\leq p <\infty$ denote by $\mathcal{J}^p$ the
Schatten $p$-class, and by $\|.\|_{\mathcal{J}^p}$ its
norm. We will use $\|.\|$ for the usual operator
norm. Using $\|A\|_{\mathcal{J}^p} = \infty$ if
$A \notin \mathcal{J}^p$, we have the following
inequalities for all operators:
\be
\|A B\|_{\mathcal{J}^p} \leq \|A\|\|B\|_{\mathcal{J}^p}, \quad
\|A B\|_{\mathcal{J}^1} \leq \|A\|_{\mathcal{J}^2} \|B|_{\mathcal{J}^2}.
\ee
Furthermore, we will use the notation of $\mathcal{J}^p$-converges to denote convergence in the
respective $\|.\|_{\mathcal{J}^p}$-norm.

The following result from \cite[Thm~IV.11.3]{tracedet} will be needed.

\begin{lemma}\label{lem:contrace}
Let $p > 0$, $A \in \mathcal{J}^p$, $T_n \xrightarrow{s} T$, $S_n \xrightarrow{s} S$ sequences
of strongly convergent bounded linear operators in some separable Hilbert space, then
\be
\|T_n A S_n^\ast - T A S^\ast\|_{\mathcal{J}^p} \rightarrow 0.
\ee
\end{lemma}

We will also need the following continuity result for $\xi$.
It will also allow us to fix the unknown constant.

\begin{lemma} \label{lem:Heps}
Suppose $H_\eps$, $\eps\in[0,1]$, is a family of self-adjoint operators, which is continuous in
the metric
\be
\rho(A,B) = \|R_A (z_0) - R_B(z_0)\|_{\mathcal{J}^1},
\ee
for some fixed $z_0\in \C\backslash\R$ and abbreviate $\xi_\eps = \xi(H_\eps, H_0)$.
Then there exists a unique choice of $\xi_\eps$ such that $\eps \mapsto \xi_\eps$ is a continuous
map $[0,1] \rightarrow L^1(\R,(\lam^2+1)^{-1}d\lam)$ with $\xi_0 = 0$.

If $H_\eps \ge \lam_0$ is bounded from below,
we can also allow $z=\lam \in (-\infty,\lam_0)$.
\end{lemma}

\begin{proof}
The first statement can be found in \cite[Lem.~8.7.5]{yafams}.
To see the second statement, let $\lam < \lam_0$ and $|\lam - z| < \lam_0 - \lam$
for some $z\in \C\backslash\R$. Abbreviate $R_\eps(z)=R_{H_\eps}(z)$.
Now using the first resolvent identity gives
\begin{align*}
\|R_\eps(z) - R_{\eps'}(z)\|_{\mathcal{J}^1} \leq &
\|R_\eps(\lam) - R_{\eps'}(\lam)\|_{\mathcal{J}^1} \\
\nn & {} + |z - \lam| \|R_\eps(z)\| \|R_\eps(\lam) - R_{\eps'}(\lam)\|_{\mathcal{J}^1}\\
\nn & {} + |z-\lam| \|R_{\eps'}(\lam)\| \|R_\eps(z) - R_{\eps'}(z)\|_{\mathcal{J}^1}
\end{align*}
and our conditions imply
$$
|z-\lam| \|R_{\eps'}(\lam)\| \leq \frac{|z-\lam|}{\lam_0 - \lam} < 1
$$
and thus
$$
\|R_\eps(z) - R_{\eps'}(z)\|_{\mathcal{J}^1} \leq
\frac{1 + \frac{|z - \lam|}{|\im(z)|}}{1 - \frac{|z-\lam|}{|\lam_0 - \lam|}}
\|R_\eps(\lam) - R_{\eps'}(\lam)\|_{\mathcal{J}^1},
$$
from which the statement follows.
\end{proof}

Our final aim is to find some conditions which allow us to verify the assumptions of this
lemma. To do this, we derive some properties of relatively bounded operators
multiplied by strongly continuous families of operators. The key example for these operators will
be multiplication operators by characteristic functions strongly 
converging to the identity operator.

\begin{hypothesis} \label{hyp:h0v}
Suppose $H_0$ and $V$ are self-adjoint such that:
\begin{itemize}
\item[(i)] $V$ is relatively bounded with respect to $H_0$ with $H_0$ bound less than one or
\item[(i')] $H_0$ is bounded from below and $|V|$ is relatively form bounded with respect
to $H_0$ with relative form bound less than one and
\item[(ii)] $|V|^{1/2} R_{H_0}(z)$ is Hilbert--Schmidt for one (and hence for all) $z\in\rho(H_0)$.
\end{itemize}
\end{hypothesis}

We recall that (i) means that $\dom(V) \supseteq \dom(H_0)$ and for some $a < 1$, $b \ge 0$,
\be
\|V \psi\| \leq a \|H_0 \psi\| + b \|\psi\|,
\quad \forall \psi \in \dom(H_0).
\ee
and (i') means that $\fdom(V) \supseteq \fdom(H_0)$ and  for some $a < 1$, $b \ge 0$,
\be \label{defrfb}
\spr{\psi}{|V| \psi} \leq a \spr{\psi}{H_0 \psi} + b \|\psi\|^2,
\quad \forall \psi \in \fdom(H_0).
\ee
Note that in case (i') all sums of operators are meant as form sums.

\begin{lemma}\label{lem:resolvconv}
Let $\eps \ni[0,1] \to K_\eps$ be a strongly
continuous family of bounded self-adjoint operators
which commute with $V$ and $0= K_0\le K_\eps \le K_1 = \id$. 

Assume Hypothesis~\ref{hyp:h0v}. Then
\be
H_\eps = H_0 + K_\eps V
\ee
are self-adjoint operators such that the assumptions of Lemma~\ref{lem:Heps}
hold.
\end{lemma}

\begin{proof}
We will abbreviate $V_\eps = K_\eps V$ and $R_\eps(z)= R_{H_\eps}(z)$.

We begin with the case where $V$ is relatively bounded with respect to $H_0$. By the Kato--Rellich Theorem
(\cite[Thm.~X.12]{reedsimon2}) $H_\eps$ is well-defined and self-adjoint. Moreover,
there is a $z$ with $\im(z) \neq 0$ such that $\| V R_0(z) \|  \le a <1$. Hence
$\| V_\eps R_0(z) \|  \le a$ and a straightforward calculation using the second resolvent identity,
$$
V R_\eps(z) = V R_0(z) (1 + V_\eps R_0(z))^{-1},
$$
shows that
$$
\| V R_\eps(z) \| \le \frac{a}{1-a}.
$$
Furthermore, again using the second resolvent identity, we have
$$
|V|^{1/2} R_\eps(z) = |V|^{1/2} R_0(z) ( 1 - V_\eps R_\eps(z)),
$$
which shows that
$$
\| |V|^{1/2} R_\eps(z) \|_{\mathcal{J}_2} \le \frac{1}{1-a} \| |V|^{1/2} R_0(z) \|_{\mathcal{J}_2}.
$$
To show $\mathcal{J}^1$-continuity at some fixed $\eps\in[0,1]$ observe
$$
R_{\eps'}(z) - R_\eps(z) = R_{\eps'}(z) |V|^{1/2} \Big( (K_{\eps'} - K_\eps)
\sgn(V) |V|^{1/2} R_\eps \Big),
$$
where the first term $R_{\eps'}(z) |V|^{1/2} \subseteq (|V|^{1/2} R_{\eps'}(z^*))^*$ is uniformly
$\mathcal{J}^2$-bounded in $\eps'$ by our previous argument and the second term
$\mathcal{J}^2$-converges to $0$ as $\eps'$ to $\eps$ by Lemma~\ref{lem:contrace}.

Now we come to the case where $V$ is relatively form bounded with respect to $H_0$.
By the KLMN theorem (\cite[Thm.~X.17]{reedsimon2}), $H_\eps$ is self-adjoint since we have
$$
|\spr{\psi}{V_\eps \psi}| = |\spr{|V|^{1/2}\psi}{K_\eps \sgn(V) |V|^{1/2}\psi}| \le \spr{\psi}{|V|\psi},
\quad \psi\in\fdom(V_\eps)= \fdom(V).
$$
Moreover, using (\ref{defrfb}) we have
$$
\| |V|^{1/2} R_\eps(-\lam)^{1/2} \|^2 \le a, \qquad \frac{b}{a}< \lam.
$$
We abbreviate:
$$
U = |V|^{1/2},\quad W =  \sgn(V) |V|^{1/2}.
$$
For $\lam>\frac{b}{a}$ we have (see \cite[Sec.~VI.3]{kato} or \cite[Thm.~II.12]{shqf})
\begin{align*}
R_\eps(-\lam) &= R_0(-\lam)^{1/2} ( 1 +C_\eps)^{-1} R_0(-\lam)^{1/2}, \\
C_\eps &= (U R_0(-\lam)^{1/2})^\ast (K_\eps W R_0(-\lam)^{1/2}).
\end{align*}
Hence, a straightforward calculation shows
\begin{align*}
R_\eps(-\lam) &= R_0(-\lam) - (U R_0(-\lam))^*
(1 +\ti{C}_\eps)^{-1} (K_\eps W R_0(-\lam)),\\
\ti{C}_\eps &= K_\eps W R_0(-\lam)^{1/2} (U R_0(-\lam)^{1/2})^\ast.
\end{align*}
By $\|\ti{C}_\eps\| \le a < 1$, one concludes that $(1 + \ti{C}_\eps)^{-1}$ exists as a bounded operator. 
Then
\begin{align*}
D_{\eps,\eps'}\psi &= (- \ti{C}_\eps(1 + \ti{C}_\eps)^{-1} -
\ti{C}_{\eps'} (1 + \ti{C}_{\eps'})^{-1})\psi \\
\nn &= (C_{\eps'} - C_\eps)(1 + C_\eps)^{-1} \psi -
C_{\eps'} D_{\eps,\eps'} \psi,
\end{align*}
where
$$
D_{\eps,\eps'} = (1 + \ti{C}_\eps)^{-1} - (1 + \ti{C}_{\eps'})^{-1}.
$$
Taking norms we obtain
$$
\|D_{\eps,\eps'}\psi\| = \frac{1}{1 - a} \|(C_{\eps'} - C_\eps)(1 + C_\eps)^{-1} \psi\|,
$$
where the last term converges to $0$ as $\eps'\to\eps$.
This implies, that $(1 + \ti{C}_\eps)^{-1}$ is strongly
continuous. Now, we obtain for the difference of resolvents
$$
R_\eps(-\lam) - R_{\eps'}(-\lam) = (U R_0(-\lam))^*
((1 +\ti{C}_\eps)^{-1} K_\eps - 
(1 +\ti{C}_{\eps'})^{-1} K_{\eps'}))
(  W R_0(-\lam))
$$
which $\mathcal{J}^1$-converges to $0$ as $\eps\to\eps'$ by
 Lemma~\ref{lem:contrace} and by $U R_0(-\lam)$ and
$W R_0(-\lam)$ being Hilbert--Schmidt.
\end{proof}

\noindent
{\bf Acknowledgments.}
We thank F.\ Gesztesy for several helpful discussions and hints with respect to the
literature. In addition, we are indebted to the referee for constructive remarks.


\begin{thebibliography}{XXXX}
\bibitem{bi} M.Sh. Birman, {\em On the spectrum of singular boundary value problems},
 AMS Translations (2) {\bf 53}, 23--80 (1966).
\bibitem{ea} M.S.P. Eastham, {\em The spectral theory of periodic differential equations},
Scottish Academic Press, Edinburgh, 1973.
\bibitem{gs} F. Gesztesy and B. Simon, {\em A short proof of Zheludev's
theorem}, Trans. Am. Math. Soc. {\bf 335}, 329--340 (1993).
\bibitem{gu} F. Gesztesy and M. \"Unal, {\em Perturbative oscillation criteria
and Hardy-type inequalities}, Math. Nachr.\ {\bf 189}, 121--144 (1998).
\bibitem{gst} F. Gesztesy, B. Simon, and G. Teschl, {\em Zeros of the Wronskian
and renormalized oscillation Theory}, Am. J. Math. {\bf 118}, 571--594 (1996).
\bibitem{tracedet} I. Gohberg, S. Goldberg, and N. Krupnik, {\em Traces and Determinants of Linear Operators},
Birkh\"auser, Basel, 2000.
\bibitem{har1} P.~Hartman, {\em Differential equations with
non-oscillatory eigenfunctions}, Duke Math.~J. {\bf 15}, 697--709 (1948).
\bibitem{har2} P.~Hartman, {\em A characterization of the spectra of
one-dimensional wave equations}, Am.~J.~Math. {\bf 71}, 915--920 (1949).
\bibitem{har3} P.~Hartman and C.R.~Putnam, {\em The least cluster
point of the spectrum of boundary value problems}, Am.~J.~Math.
{\bf 70}, 849--855 (1948).
\bibitem{kato} T. Kato, {\em Perturbation Theory for Linear Operators}, Springer, New York, 1966.
\bibitem{kn} A. Kneser, {\em Untersuchungen \"uber die reellen
Nullstellen der Integrale linearer Differentialgleichungen}, Math. Ann. {\bf 42},
409--435 (1893).
\bibitem{khr} S.V. Khryashchev, {\em Discrete spectrum for a periodic Schr\"odinger operator
perturbed by a decreasing potential}, Operator Theory: Adv. and Appl. {\bf 46}, 109--114 (1990).
\bibitem{krein} M.G. Krein, {\em Perturbation determinants and a formula for the traces of unitary and
self-adjoint operators}, Sov. Math. Dokl. {\bf 3}, 707--710 (1962).
\bibitem{kt2} H. Kr\"uger and G. Teschl, {\em Relative oscillation theory for Sturm--Liouville
operators extended}, J. Funct. Anal. {\bf 254-6}, 1702--1720 (2008).
\bibitem{kt3} H. Kr\"uger and G. Teschl, {\em Effective Pr\"ufer angles and relative oscillation criteria},
arXiv:0709.0127.
\bibitem{lei} W. Leighton, {\em On self-adjoint differential
equations of second order}, J.~London Math.~Soc. {\bf 27}, 37--47 (1952).
\bibitem{reedsimon2} M. Reed and B. Simon, {Methods of Modern Mathematical Physics II. Fourier Analysis,
Self-Adjointness}, Academic Press, New York, 1975.
\bibitem{rb1} F. S. Rofe-Beketov, {\em A test for the finiteness of the number
of discrete levels introduced into gaps of a continuous spectrum by
perturbations of a periodic potential}, Soviet Math. Dokl. {\bf 5},
689--692 (1964).
\bibitem{rb2} F.S. Rofe-Beketov, {\em Spectral analysis of the Hill operator and its perturbations},
FunkcionalÕny\"i analiz {\bf 9}, 144--155 (1977) (Russian).
\bibitem{rb3} F.S. Rofe-Beketov, {\em A generalisation of the Pr\"ufer transformation and the discrete
spectrum in gaps of  the continuous one}, Spectral Theory of Operators, 146--153, Baku, Elm, 1979 (Russian).
\bibitem{rb4} F.S. Rofe-Beketov, {\em Spectrum perturbations, the Kneser-type constants and the
effective masses of zones-type potentials}, Constructive Theory of Functions Õ84, 757--766, Sofia, 1984.
\bibitem{rb5} F.S. Rofe-Beketov, {\em Kneser constants and effective masses for band potentials},
Sov. Phys. Dokl. {\bf 29}, 391--393 (1984).
\bibitem{rbk} F.S. Rofe-Beketov and A.M. Kholkin, {\em Spectral analysis of differential operators.
Interplay between spectral and oscillatory properties}, World Scientific, Hackensack, 2005.
\bibitem{kms} K.M. Schmidt, \textit{Critical coupling constants and eigenvalue asymptotics of perturbed periodic
Sturm--Liouville operators}, Commun. Math. Phys. 211, 465--485 (2000).
\bibitem{kms2} K.M. Schmidt, \textit{Relative oscillation non-oscillation criteria for 
perturbed periodic Dirac systems}, J. Math. Anal. Appl. {\bf 246}, 591--607 (2000). 
\bibitem{kms3} K.M. Schmidt, \textit{An application of Gesztesy-Simon-Teschl oscillation theory to a
problem in differential geometry}, J. Math. Anal. Appl. {\bf 261}, 61--71 (2001). 
\bibitem{shqf} B. Simon, {\em Quantum Mechanics for Hamiltonians Defined as Quadratic Forms},
Princeton University Press, Princeton, 1971.
\bibitem{str} B. Simon, {\em Trace Ideals and Their Applications}, 2nd ed., Amer. Math. Soc., Providence, 2005.
\bibitem{sro} B. Simon, \textit{Spectral Analysis of Rank One Perturbations and Applications}, Lecture notes
from Vancouver Summer School in Mathematical Physics, August 10--14, 1993.
\bibitem{sim} B. Simon, {\em Sturm oscillation and comparison theorems}, in Sturm--Liouville Theory: Past and Present
(eds. W. Amrein, A. Hinz and D. Pearson), 29--43, Birkh\"auser, Basel, 2005.
\bibitem{sw} G. Stolz and J. Weidmann, {\em Approximation of isolated eigenvalues of
ordinary differential operators}, J. Reine und Angew. Math. {\bf 445}, 31--44 (1993).
\bibitem{stu} J.C.F. Sturm, {\em M\'emoire sur les \'equations
diff\'erentielles lin\'eaires du second ordre}, J.~Math.~Pures
Appl., {\bf 1}, 106--186 (1836).
\bibitem{tosc} G. Teschl, {\em Oscillation theory and renormalized oscillation
theory for Jacobi operators}, J. Diff. Eqs. {\bf 129}, 532--558 (1996).
\bibitem{toscd} G. Teschl, {\em Renormalized oscillation theory for Dirac operators}, Proc. Amer.
Math. Soc. {\bf 126}, 1685--1695 (1998).
\bibitem{tjac} G. Teschl, {\em Jacobi Operators and Completely Integrable Nonlinear Lattices},
Math. Surv. and Mon. {\bf 72}, Amer. Math. Soc., Rhode Island, 2000.
\bibitem{te} G. Teschl, {\em On the approximation of isolated eigenvalues of ordinary differential
operators}, Proc. Amer. Math. Soc. {\bf 136}, 2473--2476 (2008).
\bibitem{troscd} G. Teschl, {\em Relative oscillation theory for Dirac operators}, in preparation.
\bibitem{troscj} G. Teschl, {\em Relative oscillation theory for Jacobi operators}, in preparation.
\bibitem{wa} W. Walter, {\em Ordinary Differential Equations}, Springer, New York, 1998.
\bibitem{wd} J. Weidmann, {\em Zur Spektraltheorie von Sturm--Liouville--Operatoren},
Math. Z. {\bf 98}, 268--302 (1967).
\bibitem{wd1} J. Weidmann, {\em Oszillationsmethoden f\"ur Systeme gew\"ohnlicher Differentialgleichungen},
Math. Z. {\bf 119}, 349--373 (1971).
\bibitem{wdln} J. Weidmann, {\em Spectral Theory of Ordinary Differential Operators},
Lecture Notes in Mathematics, {\bf 1258}, Springer, Berlin, 1987.
\bibitem{wd2} J. Weidmann, {\em Spectral theory of Sturm--Liouville operators; approximation by
regular problems}, in Sturm--Liouville Theory: Past and Present (eds. W. Amrein, A. Hinz and
D. Pearson), 29--43, Birkh\"auser, Basel, 2005.
\bibitem{yafams} D.R. Yafaev, {\em Mathematical Scattering Theory: General Theory},
Amer. Math. Soc., Rhode Island, 1992.
\bibitem{zet} A. Zettl, {\em Sturm--Liouville Theory}, Amer. Math. Soc., Rhode Island, 2005.
\bibitem{zhe} V.A. Zheludev, {\em Perturbation of the spectrum of the one-dimensional self-adjoint
Schr\"odinger operator with a periodic potential},
Topics in Mathematical Physics, Vol. 4: M.Sh. Birman (ed), 55--76, Consultants Bureau, New York, 1971. 
\end{thebibliography}
\end{document}